\documentclass[twoside]{mystyle}

\newtheorem{theorem}{Theorem}[section]
\newtheorem{corollary}[theorem]{Corollary}

\newcommand{\qed}{\mbox{$\Box$}}
\renewcommand{\qed}{\hskip0.25em\raisebox{0.6ex}{\framebox[0.5em][l]{\ }}\vspace{.5pc}}

\newcommand{\twospace}{\renewcommand{\baselinestretch}{1.5}\normalsize}

\newcommand{\num}{\renewcommand{\theequation}{\thetheorem}\addtocounter{theorem}{1}}

\def\theequation{\thetheorem}

\NONUMBIB

\begin{document}

\twospace

\SPECFNSYMBOL{}{}{}{}{}{}{}{}{}%

\AOPMAKETITLE

\AOPAMS{Primary 60J10, 60B10; secondary 20E22.}
\AOPKeywords{Markov chain, random walk, rate of convergence to
stationarity, mixing time, wreath product, Bernoulli--Laplace
diffusion, complete monomial group, hyperoctahedral group, homogeneous
space, M\"{o}bius inversion.}
\AOPtitle{MIXING TIMES FOR MARKOV CHAINS ON WREATH PRODUCTS
AND RELATED HOMOGENEOUS SPACES}
\AOPauthor{James~Allen~Fill and Clyde~H.~Schoolfield,~Jr.}
\AOPaffil{The Johns Hopkins University and Harvard University}
\AOPlrh{J.~FILL AND C.~SCHOOLFIELD}
\AOPrrh{MIXING TIMES}

\AOPAbstract{We develop a method for analyzing the mixing times for a 
quite general class of Markov chains on the complete monomial group $G \wr
S_n$ and a quite general class of Markov chains on the homogeneous space
$(G \wr S_n) / (S_r \times S_{n-r})$.  We derive an exact formula for the
$L^2$ distance in terms of the $L^2$ distances to uniformity for closely
related random walks on the symmetric groups $S_j$ for $1 \leq j \leq n$
or for closely related Markov chains on the homogeneous spaces $S_{i+j}/
(S_i \times S_j)$ for various values of~$i$ and~$j$, respectively.  Our
results are consistent with those previously known, but our method is
considerably simpler and more general.}

\maketitle

\BACKTONORMALFOOTNOTE{3}

\thispagestyle{empty}

\twospace

\section{Introduction and Summary.} \label{sec1}

In the proofs of many of the results of Schoolfield (1999a), the~$L^2$
distance to uniformity for the random walk (on the so-called \emph{wreath
product\/} of a group~$G$ with the symmetric group~$S_n$) being analyzed
is often found to be expressible in terms of the~$L^2$ distance to
uniformity for related random walks on the symmetric groups~$S_j$ with $1
\leq j \leq n$.  Similarly, in the proofs of many of the results of
Schoolfield~(1999b), the~$L^2$ distance to stationarity for the Markov
chain being analyzed is often found to be expressible in terms of
the~$L^2$ distance to stationarity of related Markov chains on the
homogeneous spaces $S_{i+j}/ (S_i \times S_j)$ for various values of~$i$
and~$j$.  It is from this observation 
that the results of this paper have evolved.  We develop a method, with
broad applications, for bounding the rate of convergence to stationarity
for a general class of random walks and Markov chains in terms of closely
related chains on the symmetric groups and related homogeneous spaces.
Certain specialized problems of this sort were previously analyzed with
the use of group representation theory.  Our analysis is more directly
probabilistic and yields some insight into the basic structure of the
random walks and Markov chains being analyzed.

\subsection{Markov Chains on $G~\wr~S_n$.} \label{sec1.1}

We now describe one of the two basic set-ups we will be considering 
[namely, the one corresponding to the results in Schoolfield~(1999a)].  
Let~$n$ be a positive integer and let~$P$ be a probability measure defined
on a finite set~$G$ ($ = \{1, \ldots, m\}$, say).  Imagine~$n$ cards,
labeled~$1$ through~$n$ on their fronts, arranged on a table in sequential
order.  Write the number~$1$ on the back of each card.  Now repeatedly
permute the cards and rewrite the numbers on their backs, as follows.  For
each independent repetition, begin by choosing integers~$i$ and~$j$
independently according to~$P$.

If $i \neq j$, transpose the cards in positions~$i$ and~$j$.  Then,
(probabilistically) independently of the choice of~$i$ and~$j$, replace 
the numbers on the backs of the transposed cards with two numbers chosen
independently from~$G$ according to~$P$.

If $i = j$ (which occurs with probability $1 / n$), leave all cards in
their current positions.  Then, again independently of the choice of~$j$,
replace the number on the back of the card in position~$j$ by a number
chosen according to~$P$.

Our interest is in bounding the mixing time for Markov chains of the sort
we have described.  More generally, consider any probability measure,
say~$\widehat{Q}$, on the set of ordered pairs~$\hat{\pi}$ of the form
$\hat{\pi} = (\pi, J)$, where $\pi$ is a permutation of $\{1, \ldots, n\}$
and~$J$ is a subset of the set of fixed points of~$\pi$.  At each time
step, we choose such a $\hat{\pi}$ according to~$\widehat{Q}$ and then
(a)~permute the cards by multiplying the current permutation of
front-labels by~$\pi$; and (b)~replace the back-numbers of all cards whose
positions have changed, and also every card whose (necessarily unchanged)
position belongs to~$J$, by numbers chosen independently according to~$P$.

The specific transpositions example discussed above fits the more general
description, taking~$\widehat{Q}$ to be defined by 
\num \begin{equation} \label{1.1}
\begin{array}{rcll}

\widehat{Q}(e,\{j\})  &  :=  &  \displaystyle \frac{1}{n^2}  &  \mbox{for
any $j \in [n]$, with~$e$ the identity permutation}, \vspace{.5pc} \\

\widehat{Q}(\tau,\emptyset)  &  :=  &  \displaystyle \frac{2}{n^2}  &  
\mbox{for any transposition~$\tau$}, \vspace{.5pc} \\

\widehat{Q}(\hat{\pi})  &  :=  &  0  &  \mbox{otherwise}.

\end{array}
\end{equation}
\noindent
When $m = 1$, i.e.,\ when the aspect of back-number labeling is ignored,
the state space of the chain can be identified with the symmetric
group~$S_n$, and the mixing time can be bounded as in the following
classical result, which is Theorem~1 of Diaconis and Shahshahani (1981)
and was later included in Diaconis (1988) as Theorem~5 in Section~D of
Chapter~3.  The total variation norm ($\|\cdot\|_{\mbox{\rm \scriptsize
TV}}$) and the~$L^2$ norm ($\|\cdot\|_2$) will be reviewed in
Section~\ref{sec1.3}.

\begin{theorem} \label{1.2}
Let~$\nu^{*k}$ denote the distribution at time~$k$ for the random 
transpositions chain~{\rm (\ref{1.1})} when $m = 1$, and let~$U$ be the
uniform distribution on~$S_n$.  Let $k = \frac{1}{2} n \log n + cn$.  Then
there exists a universal constant $a > 0$ such that
\[
\| \nu^{*k} - U \|_{\mbox{\rm \scriptsize TV}} \ \ \leq \ \
\mbox{$\frac{1}{2}$}\ \| \nu^{*k} - U \|_2 \ \
\leq \ \ a e^{- 2 c} \ \ \ \mbox{\rm for all\ }c > 0.
\]
\end{theorem}

Without reviewing the precise details, we remark that this bound is sharp,
in that there is a matching lower bound for total variation (and hence
also for $L^2$).  Thus, roughly put, $\frac{1}{2} n \log n + cn$ steps are
necessary and sufficient for approximate stationarity.

Now consider the chain~(\ref{1.1}) for general $m \geq 2$, but restrict
attention to the case that~$P$ is uniform on~$G$.  An elementary approach
to bounding the mixing time is to combine the mixing time result of
Theorem~\ref{1.2} (which measures how quickly the cards get mixed up)
with a coupon collector's analysis (which measures how quickly their
back-numbers become random).  This approach is carried out in
Theorem~3.6.4 of Schoolfield~(1999a), but gives an upper bound only on
total variation distance.  If we are to use the chain's mixing-time
analysis in conjunction with the powerful comparison technique of Diaconis
and Saloff-Coste~(1993a, 1993b) to bound mixing times for other more
complicated chains, as is done for example in Section~4 of
Schoolfield~(1999a), we need an upper bound on~$L^2$ distance.

Such a bound can be obtained using group representation theory.  Indeed,
the Markov chain we have described is a random walk on the complete
monomial group $G \wr S_n$, which is the wreath product of the group~$G$
with~$S_n$; see Schoolfield~(1999a) for further background and discussion.
The following result is Theorem~3.1.3 of Schoolfield~(1999a).

\begin{theorem} \label{1.3}
Let~$\nu^{*k}$ denote the distribution at time~$k$ for the random 
transpositions chain~{\rm (\ref{1.1})} when~$P$ is uniform on~$G$ (with
$|G| \geq 2$).  Let $k = \frac{1}{2} n \log n + \frac{1}{4} n \log(|G| -
1) + cn$.  Then there exists a universal constant $b > 0$ such that
\[
\| \nu^{*k} - U \|_{\mbox{\rm \scriptsize TV}} \ \ \leq \ \ 
\mbox{$\frac{1}{2}$} \ \| \nu^{*k} - U \|_2 \ \ \leq \ \ b e^{- 2 c} \ \ \
\mbox{\rm for all\ }c > 0.
\]
\end{theorem}

For~$L^2$ distance (but not for TV distance), the presence of the 
additional term $\frac{1}{4} n \log(|G| - 1)$ in the mixing-time bound is
``real,'' in that there is a matching lower bound:\ see the table at the
end of Section~3.6 in Schoolfield~(1999a).

The group-representation approach becomes substantially more difficult to
carry out when the card-rearrangement scheme is something other than
random transpositions, and prohibitively so if the resulting 
step-distribution on~$S_n$ is not constant on conjugacy classes.  
Moreover, there is no possibility whatsoever of using this approach 
when~$P$ is non-uniform, since then we are no longer dealing with random
walk on a group.

In Section~\ref{sec2} we provide an $L^2$-analysis of our chain for 
completely general shuffles~$\widehat{Q}$ of the sort we have described.
More specifically, in Theorem~\ref{2.3} we derive an exact formula for
the~$L^2$ distance to stationarity in terms of the $L^2$ distance for
closely related random walks on the symmetric groups $S_j$ for $1 \leq j
\leq n$.  Subsequent corollaries establish more easily applied results in
special cases.  In particular, Corollary~\ref{2.8} extends
Theorem~\ref{1.3} to handle non-uniform $P$.

Our new method does have its limitations.  The back-number randomizations
must not depend on the current back numbers (but rather chosen afresh
from~$P$), and they must be independent and identically distributed from
card to card.  So, for example, we do not know how to adapt our method to
analyze the ``paired-shuffles'' random walk of Section~3.7 in Schoolfield
(1999a).

\subsection{Markov Chains on $(G~\wr~S_n) / (S_r~\times~S_{n - r})$.}
\label{sec1.2}

We now turn to our second basic set-up [namely, the one corresponding to
the results in Schoolfield~(1999b)].  Again, let~$n$ be a positive integer
and let~$P$ be a probability measure defined on a finite set~$G = \{1,
\ldots, m\}$.

Imagine two racks, the first with positions labeled~$1$ through~$r$ and
the second with positions labeled $r + 1$ through~$n$.  Without loss of
generality, we assume that $1 \leq r \leq n/2$.  Suppose that there
are~$n$ balls, labeled with serial numbers~$1$ through~$n$, each initially
placed at its corresponding rack position.  On each ball is written the
number~$1$, which we shall call its $G$-number.  Now repeatedly rearrange
the balls and rewrite their $G$-numbers, as follows.

Consider any~$\widehat{Q}$ as in Section~\ref{sec1.1}.  At each time step, 
choose~$\hat{\pi}$ from~$\widehat{Q}$ and then~(a) permute the balls by 
multiplying the current permutation of serial numbers by~$\pi$;
(b)~independently, replace the $G$-numbers of all balls whose positions
have changed as a result of the permutation, and also every ball whose
(necessarily unchanged) position belongs to~$J$, by numbers chosen
independently from~$P$; and (c)~rearrange the balls on each of the two
racks so that their serial numbers are in increasing order.

Notice that steps (a)--(b) are carried out in precisely the same way
as steps (a)--(b) in Section~\ref{sec1.1}.
The state of the system is completely determined, at each step, by the 
ordered $n$-tuple of $G$-numbers of the~$n$ balls $1, 2, \ldots, n$ and
the unordered set of serial numbers of balls on the first rack.  We have
thus described a Markov chain on the set of all $|G|^n \cdot {n \choose
r}$ ordered pairs of $n$-tuples of elements of $G$ and $r$-element subsets
of a set with~$n$ elements.

In our present setting, the transpositions example (\ref{1.1}) fits the
more general description, taking~$\widehat{Q}$ to be defined by 
\num \begin{equation} \label{1.4}
\begin{array}{rcll}

\widehat{Q}(\kappa,\{j\})  &  :=  &  \displaystyle \frac{1}{n^2 r! (n-r)!}
&  \mbox{where $\kappa \in K$ and $j \in [n]$}, \vspace{.5pc} \\

\widehat{Q}(\kappa,\{i,j\})  &  :=  &  \displaystyle \frac{2}{n^2 r!
(n-r)!}  &  \begin{array}{l} \mbox{where $\kappa \in K$ and $i \neq j$}
\\  \mbox{with $i,j \in [r]$ or $i,j \in [n] \setminus [r]$}, \end{array}
\vspace{.5pc} \\

\widehat{Q}(\tau \kappa,\emptyset)  &  :=  &  \displaystyle \frac{2}{n^2
r! (n-r)!}  &  \mbox{where $\tau \kappa \in T K$}, \vspace{.5pc} \\

\widehat{Q}(\hat{\pi})  &  :=  &  0  &  \mbox{otherwise},

\end{array}
\end{equation}
\noindent
where $K := S_r \times S_{n-r}$, $T$ is the set of all
transpositions in $S_n \setminus K$, and $T K := \{ \tau \kappa \in S_n :
\tau \in T \ \mbox{and}\ \kappa \in K \}$.  When $m = 1$, the state space
of the chain can be identified with the homogeneous space $S_n / (S_r
\times S_{n - r})$.  The chain is then a variant of the celebrated 
Bernoulli--Laplace diffusion model.  For the classical model, Diaconis and
Shahshahani~(1987) determined the mixing time.  Similarly, 
Schoolfield~(1999b) determined the mixing time of the present variant, 
which slows down the classical chain by a factor of $\frac{n^2}{2 r (n -
r)}$ by not forcing two balls to switch racks at each step.  The following
result is Theorem~2.5.3 of Schoolfield~(1999b).

\begin{theorem} \label{1.5}
Let~$\widetilde{\nu^{*k}}$ denote the distribution at time~$k$ for the
variant
{\rm (\ref{1.4})} 
of the Bernoulli--Laplace model when $m = 1$, and 
let~$\widetilde{U}$ be the uniform distribution on~$S_n / (S_r \times S_{n
- r})$.  Let $k = \frac{1}{4} n (\log n + c)$.  Then there exists a
universal constant $a > 0$ such that
\[
\| \widetilde{\nu^{*k}} - \widetilde{U} \|_{\mbox{\rm \scriptsize TV}}
\ \ \leq \ \ \mbox{$\frac{1}{2}$} \ \| \widetilde{\nu^{*k}}  -
\widetilde{U} \|_2 \ \ \leq \ \ a e^{- 2 c} \ \ \ \mbox{\rm for all\ }c >
0.
\]
\end{theorem}

Again there are matching lower bounds, for $r$ not too far from $n/2$,
so this Markov chain is twice as
fast to converge as the random walk of Theorem~\ref{1.2}.

The following analogue, for the special case $m = 2$, of 
Theorem~\ref{1.3} in the present setting was obtained as Theorem~3.1.3
of Schoolfield~(1999b).

\begin{theorem} \label{1.6}
Let~$\widetilde{\nu^{*k}}$ denote the distribution at time~$k$ for the
variant
{\rm (\ref{1.4})}
of the Bernoulli--Laplace model
when~$P$ is
uniform on~$G$
with $|G| = 2$.  Let $k = \frac{1}{4} n (\log n + c)$.  
Then there exists a universal constant $b > 0$ such that
\[
\| \widetilde{\nu^{*k}} - \widetilde{U} \|_{\mbox{\rm \scriptsize TV}}
\ \ \leq \ \ \mbox{$\frac{1}{2}$} \ \| \widetilde{\nu^{*k}} -
\widetilde{U} \|_2 \ \ \leq \ \ be^{-c/2} \ \ \ \mbox{\rm for all\ }\ c >
0. 
\]
\end{theorem}

Notice that Theorem~\ref{1.6} provides (essentially) the same mixing 
time bound as that found in Theorem~\ref{1.5}.  Again there are matching
lower bounds, for $r$ not too far from $n/2$, so this Markov chain
is twice as fast to converge as the random walk of
Theorem~\ref{1.3} in the special case $m = 2$.

In Section~\ref{sec3}, we provide a general $L^2$-analysis of our chain, 
which has state space equal to the homogeneous space $(G~\wr~S_n) / (S_r
\times S_{n-r})$.  More specifically, in Theorem~\ref{3.3} we derive an
exact formula for the~$L^2$ distance to stationarity in terms of the~$L^2$
distance for closely related Markov chains on the homogeneous spaces
$S_{i+j}/ (S_i \times S_j)$ for various values of~$i$ and~$j$.  Subsequent
corollaries establish more easily applied results in special cases.  In
particular, Corollary~\ref{3.8} extends Theorem~\ref{1.6} to handle
non-uniform $P$.

Again, our method does have its limitations.  For example,
we do not know how to adapt our method to analyze the ``paired-flips'' 
Markov chain of Section~3.4 in Schoolfield (1999b).

\subsection{Distances Between Probability Measures.} \label{sec1.3}

We now review several ways of measuring distances between
probability measures on a finite set~$G$.  Let~$R$ be a fixed
reference probability measure on~$G$ with $R(g) > 0$ for all $g \in G$.
As discussed in Aldous and Fill~(200x), for each $1 \leq p < \infty$
define the \emph{$L^p$ norm\/} $\| \nu \|_p$ of any signed measure~$\nu$ on~$G$
(with respect to~$R$)
by
\[ 
\| \nu \|_p \ \ := \ \ \left( \mathbb{E}_R
\left| \frac{\nu}{R} \right|^p \right)^{1/p} \ \ =
\ \ \left( \sum_{g \in G} \frac{\left|
\nu(g) \right|^p}{R(g)^{p-1}} \right)^{1/p}. 
\]
Thus the $L^p$ distance between any two probability measures~$P$ and~$Q$ on~$G$
(with respect to~$R$) is
\[ 
\| P - Q \|_p \ \ = \ \ \left( \mathbb{E}_R \left |\frac{
P - Q}{R} \right|^p \right)^{1/p} \ \ =
\ \ \left( \sum_{g \in G} \frac{\left|
P(g) - Q(g) \right|^p}{R(g)^{p-1}} \right)^{1/p} 
\]
\noindent
Notice that
\[ 
\| P - Q \|_1 \ \ = \ \ \sum_{g \in G} |P(g) - Q(g)|. 
\]
\noindent
In our applications we will always take $Q = R$ (and
$R$ will always be the stationary distribution
of the Markov chain under consideration at that time).
In that case, when $U$ is the uniform distribution on $G$,
\[ 
\| P - U \|_2 \ \ = \ \ \left( |G| \sum_{g \in G} |P(g) - U(g)|^2
\right)^{1/2}. 
\]

The \emph{total variation distance} between $P$ and $Q$ is defined by
\[ 
\| P - Q \|_{\mbox{\rm \scriptsize TV}} \ \ := \ \ \max_{A \subseteq G}
|P(A) - Q(A)|.
\]
\noindent
Notice that $\| P - Q \|_{\mbox{\rm \scriptsize TV}} \ = \ \frac{1}{2} \|
P - Q \|_1$.  It is a direct consequence of the Cauchy-Schwarz inequality that
\[ 
\| P - U \|_{\mbox{\rm \scriptsize TV}} \ \ \leq \ \ \displaystyle
\mbox{$\frac{1}{2}$} \ \| P - U \|_2. 
\]

If \textbf{P}$(\cdot,\cdot)$ is a reversible transition matrix on $G$ with
stationary distribution $R = {\bf P}^{\infty}(\cdot)$, then, for any $g_0
\in G$,
\[ 
\| {\bf P}^{k}\left( g_0, \cdot \right) - 
{\bf P}^{\infty}\left( \cdot \right) \|_2^2 \ \ =
\ \ \displaystyle \frac{ {\bf P}^{2k}\left( g_0, g_0 \right) }
{ {\bf P}^{\infty}\left( g_0 \right) } \ - \ 1. 
\]

All of the distances we have discussed here are indeed metrics on the space of
probability measures on~$G$.

\section{Markov Chains on $G~\wr~S_n$.} \label{sec2}

We now analyze a very general Markov chain on the complete monomial group
$G~\wr~S_n$.  It should be noted that, in the results which follow, there
is no essential use of the group structure of $G$.  So the results of this
section extend simply; in general, the Markov chain of interest is on the
set $G^n \times S_n$.

\subsection{A Class of Chains on $G~\wr~S_n$.} \label{sec2.1}

We introduce a generalization of permutations $\pi \in S_n$ which will
provide an extra level of generality in the results that follow.  Recall
that any permutation $\pi \in S_n$ can be written as the product of
disjoint cyclic factors, say
\[ 
\pi = (i_1^{(1)}\ i_2^{(1)}\ \cdots\ i_{k_1}^{(1)})\ (i_1^{(2)}\
i_2^{(2)}\ \cdots\ i_{k_2}^{(2)})\ \cdots\ (i_1^{(\ell)}\ i_2^{(\ell)}\
\cdots\ i_{k_\ell}^{(\ell)}), 
\]
\noindent
where the $K := k_1 + \cdots + k_{\ell}$ numbers $i_b^{(a)}$ are distinct
elements from $[n] := \{1, 2, \ldots, n \}$ and we may suppose $k_a \geq
2$ for $1 \leq a \leq \ell$.  The $n - K$ elements of $[n]$ not included
among the $i_b^{(a)}$ are each fixed by $\pi$; we denote this $(n-K)$-set
by $F(\pi)$.

We refer to the ordered pair of a permutation $\pi \in S_n$ and a subset
$J$ of $F(\pi)$ as an \emph{augmented permutation}.  We denote the set of
all such ordered pairs $\hat{\pi} = ( \pi, J )$, with $\pi \in S_n$ and $J
\subseteq F(\pi)$, by $\widehat{S}_n$.  For example, $\hat{\pi} \in
\widehat{S}_{10}$ given by $\hat{\pi} \ = \ \left( (12)(34)(567), \{8,10\}
\right)$ is the augmentation of the permutation $\pi \ = \ (12)(34)(567)
\in S_{10}$ by the subset $\{8,10\}$ of $F(\pi) = \{8,9,10\}$.  Notice
that any given $\hat{\pi} \in \widehat{S}_n$ corresponds to a unique
permutation $\pi \in S_n$; denote the mapping $\hat{\pi} \mapsto \pi$ by
$T$.  For $\hat{\pi} = (\pi,J) \in \widehat{S}_n$, define $I(\hat{\pi})$
to be the set of indices $i$ \emph{included} in $\hat{\pi}$, in the sense
that either $i$ is \emph{not} a fixed point of $\pi$ or $i \in J$; for our
example, $I(\hat{\pi}) \ = \ \{ 1,2,3,4,5,6,7,8,10 \}$.

\setcounter{theorem}{-1}

Let $\widehat{Q}$ be a probability measure on $\widehat{S}_n$ such that
\num \begin{equation} \label{2.0}
\widehat{Q}(\pi,J) \ = \ \widehat{Q}(\pi^{-1},J) \ \ \ \mbox{for all $\pi
\in S_n$ and $J \subseteq F(\pi) = F(\pi^{-1})$}.
\end{equation}
\noindent
We refer to this property as \emph{augmented symmetry}.  This terminology
is (in part) justified by the fact that if $\widehat{Q}$ is augmented
symmetric, then the measure $Q$ on $S_n$ induced by $T$ is given by
\[ 
Q(\pi) \ = \ \sum_{J \subseteq F(\pi)} \widehat{Q}\left( (\pi,J)
\right) \ = \ Q(\pi^{-1}) \ \ \ \mbox{for each $\pi \in S_n$} 
\]
\noindent
and so is symmetric in the usual sense.  We assume that $Q$ is not
concentrated on a subgroup of $G$ or a coset thereof.  Thus $Q^{*k}$
approaches the uniform distribution $U$ on $S_n$ for large $k$.  

Suppose that $G$ is a finite group.  Label the elements of $G$ as $g_1,
g_2, \ldots, g_{|G|}$.  Let $P$ be a probability measure defined on $G$.
Define $p_i := P(g_i)$ for $1 \leq i \leq |G|$.  To avoid trivialities, we
suppose $\displaystyle p_{\min} := \min \left\{ p_i : 1 \leq i \leq |G|
\right\} > 0$.

Let $\hat{\xi}_1, \hat{\xi}_2, \ldots$ be a sequence of independent
augmented permutations each distributed according to $\widehat{Q}$.  These
correspond uniquely to a sequence $\xi_1, \xi_2, \ldots$ of permutations
each distributed according to $Q$.  Define \textbf{Y} $:= (Y_0, Y_1, Y_2,
\ldots)$ to be the random walk on $S_n$ with $Y_0 := e$ and $Y_k := \xi_k
\xi_{k-1} \cdots \xi_1$ for all $k \geq 1$.  (There is no loss of 
generality in defining $Y_0 := e$, as any other $\pi \in S_n$ can be
transformed to the identity by a permutation of the labels.)

Define \textbf{X} $:= (X_0, X_1, X_2, \ldots)$ to be the Markov chain on
$G^n$ such that $X_0 := \vec{x}_0 = (\chi_1, \ldots, \chi_n)$ with
$\chi_i \in G$ for $1 \leq i \leq n$ and, at each step $k$ for $k \geq 1$,
the entries of $X_{k-1}$ whose positions are included in $I(\hat{\xi}_k)$
are independently changed to an element of $G$ distributed according to
$P$.  

Define \textbf{W} $:= (W_0, W_1, W_2, \ldots)$ to be the Markov chain on
$G~\wr~S_n$ such that $W_k := (X_k; Y_k)$ for all $k \geq 0$.  Notice that
the random walk on $G~\wr~S_n$ analyzed in Theorem~\ref{1.3} is a
special case of \textbf{W}, with $P$ being the uniform distribution and
$\widehat{Q}$ being defined as at (\ref{1.1}).  Let \textbf{P}$(\cdot,
\cdot)$ be the transition matrix for \textbf{W} and let
\textbf{P}$^{\infty}(\cdot)$ be the stationary distribution for
\textbf{W}.  

Notice that 
\[ 
\mbox{\textbf{P}}^{\infty} \left( \vec{x};\pi \right) \ = \
\frac{1}{n!} \prod_{i=1}^n p_{x_i} 
\]
\noindent
for any $(\vec{x};\pi) \in G~\wr~S_n$ and that 
\[ 
\mbox{\textbf{P}} \left( (\vec{x};\pi), (\vec{y};\sigma) \right) \ = \
\sum_{\hat{\rho} \in \widehat{S}_n : T(\hat{\rho}) = \sigma \pi^{-1}}
\widehat{Q}(\hat{\rho}) \left[ \prod_{j \in I(\hat{\rho})} p_{y_j} \right]
\cdot \left[ \prod_{\ell \not\in I(\hat{\rho})} \mathbb{I}(x_{\ell} =
y_{\ell}) \right] 
\]
\noindent 
for any $(\vec{x};\pi), (\vec{y};\sigma) \in G~\wr~S_n$.  Thus, using the
augmented symmetry of $\widehat{Q}$,
\[ \begin{array}{l}

\displaystyle \mbox{\textbf{P}}^{\infty}\left( \vec{x};\pi \right)
\mbox{\textbf{P}} \left( (\vec{x};\pi), (\vec{y};\sigma) \right)
\vspace{1pc} \\

\ \ \ = \displaystyle \left[ \frac{1}{n!} \prod_{i=1}^n p_{x_i} \right]
\sum_{\hat{\rho} \in \widehat{S}_n : T(\hat{\rho}) = \sigma \pi^{-1}}
\widehat{Q}(\hat{\rho}) \left[ \prod_{j \in I(\hat{\rho})} p_{y_j} \right]
\cdot \left[ \prod_{\ell \not\in I(\hat{\rho})} \mathbb{I}(x_{\ell} =
y_{\ell}) \right] \vspace{1pc} \\

\ \ \ = \displaystyle \sum_{\hat{\rho} \in \widehat{S}_n : T(\hat{\rho}) =
\sigma \pi^{-1}} \widehat{Q}(\hat{\rho}) \left[ \frac{1}{n!} \left(
\prod_{i \in I(\hat{\rho})} p_{x_i} \right) \left( \prod_{i \not\in
I(\hat{\rho})} p_{x_i} \right) \right] \cdot \left[ \prod_{j \in
I(\hat{\rho})} p_{y_j} \right] \cdot \left[ \prod_{\ell \not\in
I(\hat{\rho})} \mathbb{I}(x_{\ell} = y_{\ell}) \right] \vspace{1pc} \\

\ \ \ = \displaystyle \sum_{\hat{\rho} \in \widehat{S}_n : T(\hat{\rho}) =
\pi \sigma^{-1}} \widehat{Q}(\hat{\rho}) \left[ \frac{1}{n!} \left(
\prod_{i \in I(\hat{\rho})} p_{x_i} \right) \left( \prod_{j \not\in
I(\hat{\rho})} p_{y_j} \right) \right] \cdot \left[ \prod_{j \in
I(\hat{\rho})} p_{y_j} \right] \cdot \left[ \prod_{\ell \not\in
I(\hat{\rho})} \mathbb{I}(y_{\ell} = x_{\ell}) \right] \vspace{1pc} \\

\ \ \ = \displaystyle \left[ \frac{1}{n!} \prod_{j=1}^n p_{y_j} \right]
\sum_{\hat{\rho} \in \widehat{S}_n : T(\hat{\rho}) = \pi \sigma^{-1}}
\widehat{Q}(\hat{\rho}) \left[ \prod_{i \in I(\hat{\rho})} p_{x_i} \right]
\cdot \left[ \prod_{\ell \not\in I(\hat{\rho})} \mathbb{I}(y_{\ell} =
x_{\ell}) \right] \vspace{1pc} \\

\ \ \ = \displaystyle \mbox{\textbf{P}}^{\infty}\left( \vec{y};\sigma
\right) \mbox{\textbf{P}} \left( (\vec{y};\sigma), (\vec{x};\pi) \right).

\end{array} \]
\noindent
Therefore, \textbf{P} is
reversible, which is a necessary condition in
order to apply the comparison technique of Diaconis and Saloff-Coste (1993a).

\subsection{Convergence to Stationarity:\ Main Result.} \label{sec2.2}

For notational purposes, let 
\num \begin{equation} \label{2.1} \mu_n(J) \ \ := \ \ \widehat{Q}\{ 
\hat{\sigma} \in \widehat{S}_n : I(\hat{\sigma}) \subseteq J\}.
\end{equation}
\noindent
For any $J \subseteq [n]$, let $S_{(J)}$ be the subgroup of $S_n$ 
consisting of those $\sigma \in S_n$ with $[n] \setminus F(\sigma) 
\subseteq J$.  If $\hat{\pi} \in \widehat{S}_n$ is random with
distribution $\widehat{Q}$,
then, when the conditioning event
$$
E := \{ I(\hat{\pi}) \subseteq J \} 
\bigg[ = \{ [n] \setminus F\left( T(\hat{\pi}) \right) \subseteq J \} \bigg]
$$
has positive probability, the probability measure induced by $T$ from the
conditional distribution (call it $\widehat{Q}_{S_{(J)}}$) of $\hat{\pi}$
given~$E$ is concentrated on
$S_{(J)}$.  Call this induced measure $Q_{S_{(J)}}$.  Notice that
$\widehat{Q}_{S_{(J)}}$, like $\widehat{Q}$, is augmented symmetric and
hence that $Q_{S_{(J)}}$ is symmetric on $S_{(J)}$.  Let $U_{S_{(J)}}$ be
the uniform measure on $S_{(J)}$.  For notational purposes, let 
\num \begin{equation} \label{2.2}
d_k(J) \ \ := \ \ |J|! \| Q_{S_{(J)}}^{*k} - U_{S_{(J)}} \|_2^2.
\end{equation}
\bigskip

\emph{Example.\/}\ Let $\widehat{Q}$ be defined as at (\ref{1.1}).
Then $\widehat{Q}$ satisfies the augmented symmetry property (\ref{2.0}).
In Corollary~\ref{2.8} we will be using~$\widehat{Q}$ to define a random walk on
$G~\wr~S_n$ which is precisely the random walk analyzed in
Theorem~\ref{1.3}.

For now, however, we will be satisfied to determine~$\widehat{Q}_{S_{(J)}}$
and~$Q_{S_{(J)}}$, where $J \subseteq [n]$.
It is easy to verify that
\[ \begin{array}{rcll}

\widehat{Q}_{S_{(J)}}(e,\{j\})  &  :=  &  \displaystyle \frac{1}{|J|^2}  &
\mbox{for each $j \in J$},  \vspace{.5pc} \\

\widehat{Q}_{S_{(J)}}((p\ q),\emptyset)  &  :=  &  \displaystyle
\frac{2}{|J|^2}  &  \mbox{for each transposition $\tau \in S_n$ with
$\{p,q\} \subseteq J$},  \vspace{.5pc} \\

\widehat{Q}_{S_{(J)}}(\hat{\pi})  &  :=  &  0  &  \mbox{otherwise},

\end{array} \]
\noindent
and hence that $\widehat{Q}_{S_{(J)}}$ is the probability measure defined at
(\ref{1.1}), but with $[n]$ changed to $J$.  Thus, roughly put, the
random walk analyzed in Theorem~\ref{1.3}, conditionally restricted to
the indices in $J$, gives a random walk ``as if $J$ were the only
indices.''
\bigskip

The following result establishes an upper bound on the total variation
distance by deriving an exact formula for $\| \mbox{\textbf{P}}^k \left(
(\vec{x}_0,e), \cdot \right) - \mbox{\textbf{P}}^{\infty} (\cdot) \|_2^2$.  

\begin{theorem} \label{2.3}
Let \emph{\textbf{W}} be the Markov chain on the complete monomial group
$G~\wr~S_n$ defined in Section~{\rm \ref{sec2.1}}.  Then

\[ \begin{array}{rcl}

\displaystyle \| \emph{\textbf{P}}^{k}\left( (\vec{x}_0;e), \cdot \right)
- \emph{\textbf{P}}^{\infty}\left( \cdot \right) \|_{\mbox{\rm \scriptsize
TV}}^2  &  \leq  &  \displaystyle \mbox{$\frac{1}{4}$} \ \|
\emph{\textbf{P}}^k \left( (\vec{x}_0,e), \cdot \right) -
\emph{\textbf{P}}^{\infty} (\cdot) \|_2^2 \vspace{1pc} \\

&  =  &  \displaystyle \mbox{$\frac{1}{4}$} \sum_{J: J \subseteq [n]}
\frac{n!}{|J|!} \left[ \prod_{i \not\in J} \left( 
\mbox{$\frac{1}{p_{\chi_i}}$} - 1 \right) \right] \mu_n(J)^{2k} \ d_k(J)
\vspace{1pc} \\

&     &  \displaystyle + \ \ \mbox{$\frac{1}{4}$} \sum_{J: J \subsetneq
[n]} \frac{n!}{|J|!} \left[ \prod_{i \not\in J} \left(
\mbox{$\frac{1}{p_{\chi_i}}$} - 1 \right) \right] \mu_n(J)^{2k}.

\end{array} \]
\noindent
where $\mu_n(J)$ and $d_k(J)$ are defined at~{\rm (\ref{2.1})} and~{\rm (\ref{2.2})},
respectively.

\end{theorem}

Before proceeding to the proof, we note the following.  In the present
setting, the argument used to prove Theorem 3.6.4 of Schoolfield (1999a)
gives the upper bound
\[ 
\| \mbox{\textbf{P}}^{k}\left( (\vec{x}_0;e), \cdot \right) -
\mbox{\textbf{P}}^{\infty}\left( \cdot \right) \|_{\mbox{\rm \scriptsize
TV}} \ \ \leq \ \ \| Q^{*k} - U_{S_n} \|_{\mbox{\rm \scriptsize TV}} \ + \
\mathbb{P} \left( T > k \right), 
\]
\noindent
where $T := \inf \left\{ k \geq 1 : H_k = [n] \right\}$ and $H_k$ is
defined as at the outset of that theorem's proof.  Theorem~\ref{2.3}
provides a similar type of upper bound, but (a) we work with $L^2$
distance instead of total variation distance and (b) the analysis is more
intricate, involving the need to consider how many steps are needed to
escape sets $J$ of positions and also the need to know $L^2$ for random
walks on subsets of $[n]$.  However, Theorem~\ref{2.3} does derive an
\emph{exact} formula for $L^2$.

\proof{Proof} For each $k \geq 1$, let $\displaystyle H_k := \bigcup_{\ell
= 1}^k I(\hat{\xi}_{\ell}) \subseteq [n]$; so $H_k$ is the (random) set of
indices included in at least one of the augmented permutations 
$\hat{\xi}_1, \ldots, \hat{\xi}_k$.  For any given $w = (\vec{x};\pi) \in
G~\wr~S_n$, let $A \subseteq [n]$ be the set of indices such that $x_i
\neq \chi_i$, where $x_i$ is the $i$th entry of $\vec{x}$ and $\chi_i$ is
the $i$th entry of $\vec{x}_0$, and let $B = [n] \setminus F(\pi)$ be the
set of indices deranged by $\pi$.  Notice that $H_k \supseteq A \cup B$.
Then
\[ \begin{array}{rcl}

\displaystyle \mathbb{P}\left( W_k = (\vec{x};\pi) \right)  &  =  &
\displaystyle \sum_{C: A \cup B \subseteq C \subseteq [n]}
\mathbb{P}\left( H_k = C, W_k = (\vec{x};\pi) \right) \vspace{1pc} \\

&  =  &  \displaystyle \sum_{C: A \cup B \subseteq C \subseteq [n]}
\mathbb{P}\left( H_k = C, Y_k = \pi \right) \cdot \mathbb{P}\left( X_k =
\vec{x} \ |\ H_k = C \right) \vspace{1pc} \\

&  =  &  \displaystyle \sum_{C: A \cup B \subseteq C \subseteq [n]}
\mathbb{P}\left( H_k = C, Y_k = \pi \right) \ \prod_{i \in C} p_{x_i}.

\end{array} \]

For any $J \subseteq [n]$, we have $\mathbb{P}\left( H_k \subseteq J, Y_k
= \pi \right) = 0$ unless $B \subseteq J \subseteq [n]$, in which case
\[ \begin{array}{rcl}

\displaystyle \mathbb{P}\left( H_k \subseteq J, Y_k = \pi \right)  &  =  &
\displaystyle \mathbb{P}\left( H_k \subseteq J \right) \ \mathbb{P}\left(
Y_k = \pi \ |\ H_k \subseteq J \right) \vspace{1pc} \\

&  =  &  \displaystyle \left( \widehat{Q}\{ \hat{\sigma} \in \widehat{S}_n
: I(\hat{\sigma}) \subseteq J\} \right) ^k \ \mathbb{P}\left( Y_k = \pi \
|\ H_k \subseteq J \right) \vspace{1pc} \\

&  =  &  \displaystyle \mu_n(J)^k \ \mathbb{P}\left( Y_k = \pi \ |\ H_k
\subseteq J \right).

\end{array} \]
\noindent
Then, by M\"{o}bius inversion [see, e.g., Stanley (1986), Section 3.7],
for any $C \subseteq [n]$ we have
\[ \begin{array}{rcl}

\displaystyle \mathbb{P}\left( H_k = C, Y_k = \pi \right)  &  =  &
\displaystyle \sum_{J: J \subseteq C} \left( -1 \right)^{|C|-|J|} \
\mathbb{P}\left( H_k \subseteq J, Y_k = \pi \right) \vspace{1pc} \\

&  =  &  \displaystyle \sum_{J: B \subseteq J \subseteq C} \left( -1
\right)^{|C|-|J|} \ \mu_n(J)^k \ \mathbb{P}\left( Y_k = \pi \ |\ H_k
\subseteq J \right).

\end{array} \]

Combining these results gives
\[ \begin{array}{rcl}

\displaystyle \mathbb{P}\left( W_k = (\vec{x};\pi) \right) &  =  &
\displaystyle \sum_{C: A \cup B \subseteq C \subseteq [n]} \sum_{J: B
\subseteq J \subseteq C} \left( -1 \right)^{|C|-|J|} \ \mu_n(J)^k \
\mathbb{P}\left( Y_k = \pi \ |\ H_k \subseteq J \right) \ \prod_{i \in C}
p_{x_i} \vspace{1pc} \\

&  =  &  \displaystyle \sum_{J: B \subseteq J \subseteq [n]} \left( -1
\right)^{|J|} \mu_n(J)^k \ \mathbb{P}\left( Y_k = \pi \ |\ H_k \subseteq J
\right) \sum_{C: A \cup J \subseteq C \subseteq [n]} \ \prod_{i \in C}
(-p_{x_i}). \vspace{1pc} \\ 

\end{array} \]
\noindent
But for any $D \subseteq [n]$, we have
\[ \begin{array}{rcl}

\displaystyle \sum_{C: D \subseteq C \subseteq [n]} \prod_{i \in C}
(-p_{x_i}) &  =  &  \displaystyle \left[ \prod_{i \in D} (-p_{x_i})
\right] \sum_{E: E \subseteq [n] \setminus D} \ \ \prod_{i \in E}
(-p_{x_i}) \vspace{1pc} \\

&  =  &  \displaystyle \left[ \prod_{i \in D} (-p_{x_i}) \right] \prod_{i
\in [n] \setminus D} (1 - p_{x_i}) \vspace{1pc} \\

&  =  &  \displaystyle \prod_{i \in [n]} \left[ 1 - \mathbb{I}_D(i) -
p_{x_i} \right]

\end{array} \]
\noindent
where (as usual) $\mathbb{I}_D(i) = 1$ if $i \in D$ and $\mathbb{I}_D(i) =
0$ if $i \not\in D$.  Therefore
\[ 
\displaystyle \mathbb{P}\left( W_k = (\vec{x};\pi) \right) \ = \
\sum_{J: B \subseteq J \subseteq [n]} \left( -1 \right)^{|J|} \mu_n(J)^{k}
\ \mathbb{P}\left( Y_k = \pi \ |\ H_k \subseteq J \right) \prod_{i=1}^n
\left[ 1 - \mathbb{I}_{A \cup J}(i) - p_{x_i} \right]. 
\]

In particular, when $(\vec{x};\pi) = (\vec{x}_0; e)$, we have $A =
\emptyset = B$ and
\[ \begin{array}{rcl}

\displaystyle \mathbb{P}\left( W_k = (\vec{x}_0;e) \right)  &  =  &
\displaystyle \sum_{J: J \subseteq [n]} \left( -1 \right)^{|J|} \mu_n(J)^k
\ \mathbb{P}\left( Y_k = e \ |\ H_k \subseteq J \right) \prod_{i=1}^n
\left[ 1 - \mathbb{I}_J(i) - p_{\chi_i} \right] \vspace{1pc} \\

&  =  &  \displaystyle \left[ \prod_{i=1}^n p_{\chi_i} \right] \sum_{J: J
\subseteq [n]} \mu_n(J)^k \ \mathbb{P}\left( Y_k = e \ |\ H_k \subseteq J
\right) \prod_{i \not\in J} \left( \mbox{$\frac{1}{p_{\chi_i}}$} - 1
\right). 

\end{array} \]
\noindent
Notice that $\displaystyle \{ H_k \subseteq J \} = \bigcap_{\ell=1}^k
\left\{ I(\hat{\xi}_{\ell}) \subseteq J \right\}$ for any $k$ and $J$.  So
$\mathcal{L}$ $\left( (Y_0, Y_1, \ldots, Y_k \ | \ H_k \subseteq J)
\right)$ is the law of a random walk on $S_n$ (through step $k$) with step
distribution $Q_{S_{(J)}}$.  Thus, using the reversibility of \textbf{P}
and the symmetry of $Q_{S_{(J)}}$,
\[ \begin{array}{l}

\| \mbox{\textbf{P}}^k \left( (\vec{x}_0,e), \cdot \right) -
\mbox{\textbf{P}}^{\infty} (\cdot) \|_2^2 \ \ = \ \ \displaystyle
\frac{n!}{\prod_{i=1}^n p_{\chi_i}} \mbox{\textbf{P}}^{2k}\left(
(\vec{x}_0;e), (\vec{x}_0;e) \right) \ - \ 1 \vspace{1pc} \\

\ \ \ = \displaystyle \ n! \ \sum_{J: J \subseteq [n]} \ \left[ \prod_{i
\not\in J} \left( \mbox{$\frac{1}{p_{\chi_i}}$} - 1 \right) \right] \
\mu_n(J)^{2k} \ \mathbb{P}\left( Y_{2k} = e \ |\ H_{2k} \subseteq J
\right) \ \ - \ \ 1 \vspace{1pc} \\

\ \ \ = \displaystyle \ n! \sum_{J: J \subseteq [n]} \ \left[ \prod_{i
\not\in J} \left( \mbox{$\frac{1}{p_{\chi_i}}$} - 1 \right) \right] \
\mu_n(J)^{2k} \ \left( \| Q_{S_{(J)}}^{*k} - U_{S_{(J)}} \|_2^2 +
\frac{1}{|J|!} \right) \ \ - \ \ 1 \vspace{1pc} \\

\ \ \ = \displaystyle \ n! \sum_{J: J \subseteq [n]} \ \left[ \prod_{i
\not\in J} \left( \mbox{$\frac{1}{p_{\chi_i}}$} - 1 \right) \right] \
\mu_n(J)^{2k} \ \frac{1}{|J|!} \left( d_k(J) + 1 \right) \ \ - \ \ 1
\vspace{1pc} \\

\ \ \ = \displaystyle \ \sum_{J: J \subseteq [n]} \frac{n!}{|J|!} \left[
\prod_{i \not\in J} \left( \mbox{$\frac{1}{p_{\chi_i}}$} - 1 \right)
\right] \mu_n(J)^{2k} \ d_k(J) \vspace{1pc} \\

\ \ \ \ \ \ \ \ + \displaystyle \ \sum_{J: J \subsetneq [n]}
\frac{n!}{|J|!} \left[ \prod_{i \not\in J} \left( 
\mbox{$\frac{1}{p_{\chi_i}}$} - 1 \right) \right] \mu_n(J)^{2k},

\end{array} \]
\noindent
from which the desired result follows. $\qed$

\subsection{Corollaries.} \label{sec2.3}

We now establish several corollaries to our main result.

\begin{corollary} \label{2.4}
Let \emph{\textbf{W}} be the Markov chain on the complete monomial group
$G~\wr~S_n$ as in Theorem~{\rm \ref{2.3}}.  For $0 \leq j \leq n$, let
\[ 
M_n(j) := \max \left\{ \mu_n(J) : |J| = j \right\} \ \ \ and \ \ \
D_k(j) := \max \left\{ d_k(J) : |J| = j \right\}. 
\]
\noindent
Also let
\[ 
B(n,k) := \max \left\{ D_k(j) : 0 \leq j \leq n \right\} = \max \left\{
d_k(J) : J \subseteq [n] \right\}. 
\]
\noindent
Then
\[ \begin{array}{rcl}

\displaystyle \| \emph{\textbf{P}}^{k}\left( (\vec{x}_0;e), \cdot \right)
- \emph{\textbf{P}}^{\infty}\left( \cdot \right) \|_{\mbox{\rm \scriptsize
TV}}^2  &  \leq  &  \displaystyle \mbox{$\frac{1}{4}$} \ \|
\emph{\textbf{P}}^k \left( (\vec{x}_0,e), \cdot \right) -
\emph{\textbf{P}}^{\infty} (\cdot) \|_2^2 \vspace{1pc} \\

&  \leq  &  \displaystyle \mbox{$\frac{1}{4}$} \ B(n,k) \sum_{j=0}^n {n
\choose j} \frac{n!}{j!} \left( \mbox{$\frac{1}{p_{\min}}$} - 1
\right)^{n-j} M_n(j)^{2k} \vspace{1pc} \\

&  &  \displaystyle + \ \ \mbox{$\frac{1}{4}$} \sum_{j=0}^{n-1} {n \choose
j} \frac{n!}{j!} \left( \mbox{$\frac{1}{p_{\min}}$} - 1 \right)^{n-j}
M_n(j)^{2k}.

\end{array} \]
\end{corollary}

\proof{Proof} Notice that 
\[ 
\prod_{i \not\in J} \left( \mbox{$\frac{1}{p_{\chi_i}}$} - 1 \right) \
\leq \ \left( \mbox{$\frac{1}{p_{\min}}$} - 1 \right)^{n - |J|}. 
\]
\noindent
The result then follows readily from Theorem~\ref{2.3}. \qed

\begin{corollary} \label{2.5}
In addition to the assumptions of Theorem~{\rm \ref{2.3}} and
Corollary~{\rm \ref{2.4}}, suppose that there exists $m > 0$ such that
$\displaystyle M_n(j) \leq \left( j/n \right)^m$ for \emph{all} $0 \leq j
\leq n$.  Let $k \geq \frac{1}{m} n \log n + \frac{1}{2m} n \log \left(
\mbox{$\frac{1}{p_{\min}}$} - 1 \right) + \frac{1}{m} cn$.  Then
\[ 
\| \emph{\textbf{P}}^{k}\left( (\vec{x}_0;e), \cdot \right) -
\emph{\textbf{P}}^{\infty}\left( \cdot \right) \|_{\mbox{\rm \scriptsize
TV}} \ \ \leq \ \ \displaystyle \mbox{$\frac{1}{2}$} \
\|\emph{\textbf{P}}^k \left( (\vec{x}_0,e), \cdot \right) -
\emph{\textbf{P}}^{\infty} (\cdot) \|_2 \ \ \leq \ \ \left( 
B\left(n,k\right) \ + \ e^{-2c} \right)^{1/2}. 
\]
\end{corollary}

\proof{Proof} It follows from Corollary~\ref{2.4} that
\num \begin{equation} \label{2.6}
\begin{array}{l}

\displaystyle \| \mbox{\textbf{P}}^{k}\left( (\vec{x}_0;e), \cdot \right)
- \mbox{\textbf{P}}^{\infty}\left( \cdot \right) \|_{\mbox{\rm \scriptsize
TV}}^2 \ \ \leq \ \ \displaystyle \mbox{$\frac{1}{4}$} \ \|
\mbox{\textbf{P}}^k \left( (\vec{x}_0,e), \cdot \right) -
\mbox{\textbf{P}}^{\infty} (\cdot) \|_2^2 \nonumber \vspace{1pc} \\

\ \ \ \leq \displaystyle \ \mbox{$\frac{1}{4}$} \ B\left(n,k\right)
\sum_{j=0}^n {n \choose j} \frac{n!}{j!} \left( 
\mbox{$\frac{1}{p_{\min}}$} - 1 \right)^{n-j} \left( \frac{j}{n}
\right)^{2km} \nonumber \vspace{1pc} \\

\ \ \ \ \ \ \ \ + \displaystyle \ \mbox{$\frac{1}{4}$} \sum_{j=0}^{n-1} {n
\choose j} \frac{n!}{j!} \left( \mbox{$\frac{1}{p_{\min}}$} - 1
\right)^{n-j} \left( \frac{j}{n} \right)^{2km}.

\end{array}
\end{equation}
\noindent
If we let $i = n-j$, then the upper bound becomes
\[ \begin{array}{l}

\displaystyle \| \mbox{\textbf{P}}^{k}\left( (\vec{x}_0;e), \cdot \right)
- \mbox{\textbf{P}}^{\infty}\left( \cdot \right) \|_{\mbox{\rm \scriptsize
TV}}^2 \ \ \leq \ \ \displaystyle \mbox{$\frac{1}{4}$} \ \|
\mbox{\textbf{P}}^k \left( (\vec{x}_0,e), \cdot \right) -
\mbox{\textbf{P}}^{\infty} (\cdot) \|_2^2 \vspace{1pc} \\

\ \ \ \leq \displaystyle \ \mbox{$\frac{1}{4}$} \ B\left(n,k\right)
\sum_{i=0}^n {n \choose i} \frac{n!}{(n-i)!} \left( 
\mbox{$\frac{1}{p_{\min}}$} - 1 \right)^i \left( 1 - \mbox{$\frac{i}{n}$}
\right)^{2km} \vspace{1pc} \\

\ \ \ \ \ \ \ \ + \displaystyle \ \mbox{$\frac{1}{4}$} \sum_{i=1}^n {n
\choose i} \frac{n!}{(n-i)!} \left( \mbox{$\frac{1}{p_{\min}}$} - 1
\right)^i \left( 1 - \mbox{$\frac{i}{n}$} \right)^{2km} \vspace{1pc} \\

\ \ \ \leq  \displaystyle \ \mbox{$\frac{1}{4}$} \ B\left(n,k\right)
\sum_{i=0}^n \frac{1}{i!} n^{2i} \left( \mbox{$\frac{1}{p_{\min}}$} - 1
\right)^i e^{-2ikm/n} \ \ + \ \ \mbox{$\frac{1}{4}$} \sum_{i=1}^n
\frac{1}{i!} n^{2i} \left( \mbox{$\frac{1}{p_{\min}}$} - 1 \right)^i
e^{-2ikm/n}. 

\end{array} \]
\noindent
Notice that if $k \geq \frac{1}{m} n \log n + \frac{1}{2m} n \log \left(
\mbox{$\frac{1}{p_{\min}}$} - 1 \right) + \frac{1}{m} cn$, then
\[ 
e^{-2ikm/n} \ \leq \ \left[ \frac{e^{-2c}}{\left(
\mbox{$\frac{1}{p_{\min}}$} - 1 \right) n^2} \right]^i, 
\]
\noindent
from which it follows that
\[ \begin{array}{rcl}

\displaystyle \| \mbox{\textbf{P}}^{k}\left( (\vec{x}_0;e), \cdot \right)
- \mbox{\textbf{P}}^{\infty}\left( \cdot \right) \|_{\mbox{\rm \scriptsize
TV}}^2  &  \leq  &  \displaystyle \mbox{$\frac{1}{4}$} \ \|
\mbox{\textbf{P}}^k \left( (\vec{x}_0,e), \cdot \right) -
\mbox{\textbf{P}}^{\infty} (\cdot) \|_2^2 \vspace{1pc} \\

&  \leq  &  \displaystyle \mbox{$\frac{1}{4}$} \ B\left(n,k\right)
\sum_{i=0}^n \frac{1}{i!} \left( e^{-2c} \right)^i \ \ + \ \
\mbox{$\frac{1}{4}$} \sum_{i=1}^n \frac{1}{i!} \left( e^{-2c} \right)^i
\vspace{1pc} \\

&  \leq  &  \displaystyle \mbox{$\frac{1}{4}$} \ B\left(n,k\right)
\exp\left( e^{-2c} \right) \ \ + \ \ \mbox{$\frac{1}{4}$} e^{-2c}
\exp\left( e^{-2c} \right).

\end{array} \]
\noindent
Since $c > 0$, we have $\exp\left(e^{-2c}\right) < e$.  Therefore
\[ 
\displaystyle \| \mbox{\textbf{P}}^{k}\left( (\vec{x}_0;e), \cdot
\right) - \mbox{\textbf{P}}^{\infty}\left( \cdot \right) \|_{\mbox{\rm
\scriptsize TV}}^2 \ \ \leq \ \ \displaystyle \mbox{$\frac{1}{4}$} \ \|
\mbox{\textbf{P}}^k \left( (\vec{x}_0,e), \cdot \right) -
\mbox{\textbf{P}}^{\infty} (\cdot) \|_2^2 \ \ \leq \ \ B\left(n,k\right) \
+ \ e^{-2c}, 
\]
\noindent
from which the desired result follows.  $\qed$

\begin{corollary} \label{2.7}
In addition to the assumptions of Theorem~{\rm \ref{2.3}} and
Corollary~{\rm \ref{2.4}}, suppose that a set with the distribution of 
$I(\hat{\sigma})$ when $\hat{\sigma}$ has distribution $\widehat{Q}$ can
be constructed by first choosing a set size $0 < \ell \leq n$ according to
a probability mass function $f_n( \cdot )$ and then choosing a set $L$
with $|L| = \ell$ \emph{uniformly} among all such choices.  Let $k \geq n
\log n + \frac{1}{2} n \log \left( \mbox{$\frac{1}{p_{\min}}$} - 1 \right)
+ cn$.  Then 
\[ 
\| \emph{\textbf{P}}^{k}\left( (\vec{x}_0;e), \cdot \right) -
\emph{\textbf{P}}^{\infty}\left( \cdot \right) \|_{\mbox{\rm \scriptsize
TV}} \ \ \leq \ \ \displaystyle \mbox{$\frac{1}{2}$} \
\|\emph{\textbf{P}}^k \left( (\vec{x}_0,e), \cdot \right) -
\emph{\textbf{P}}^{\infty} (\cdot) \|_2 \ \ \leq \ \ \left( B(n,k) \ + \
e^{-2c} \right)^{1/2}. 
\]
\end{corollary}

\proof{Proof} We apply Corollary~\ref{2.5}.  Notice that
\[ 
\widehat{Q}\{ \hat{\sigma} \in \widehat{S}_n : I(\hat{\sigma}) = L\} \
= \ \left\{ \begin{array}{ll} f_n(\ell)/ {n \choose \ell}  &  \mbox{if
$|L| = \ell$},  \\  0  &  \mbox{otherwise}. \end{array} \right.  
\]
\noindent
Then, for any $J \subseteq [n]$ with $|J| = j$,
\[ \begin{array}{rcl}

M_n(j) &  =  &  \displaystyle \widehat{Q}\{ \hat{\sigma} \in \widehat{S}_n
: I(\hat{\sigma}) \subseteq J\} \ \ = \ \ \sum_{L \subseteq J}
\widehat{Q}\{ \hat{\sigma} \in \widehat{S}_n : I(\hat{\sigma}) = L\}
\vspace{1pc} \\

&  =  &  \displaystyle \sum_{\ell=1}^j \frac{{j \choose \ell}
f_n(\ell)}{{n \choose \ell}} \ \ \leq \ \ \frac{j}{n} \sum_{\ell=1}^j 
f_n(\ell) \ \ \leq \ \ \frac{j}{n}.

\end{array} \]
\noindent
The result thus follows from Corollary~\ref{2.5}, with $m = 1$. \qed

Theorem~\ref{2.3}, and its subsequent corollaries, can be used to bound
the distance to stationarity of many different Markov chains \textbf{W} on
$G~\wr~S_n$ for which bounds on the $L^2$ distance to uniformity for
the related random walks on $S_j$ for $1 \leq j \leq n$ are known.
Theorem~\ref{1.2} provides such bounds for random walks generated by
random transpositions, showing that $\frac{1}{2} j \log j$ steps are
sufficient.  Roussel (1999) has studied random walks on $S_n$ generated
by permutations with $n-m$ fixed points for $m = 3, 4, 5, \mathrm{and\ }
6$.  She has shown that $\frac{1}{m} n \log n$ steps are both necessary
and sufficient.

Using Theorem~\ref{1.2}, the following result establishes an upper
bound on both the total variation distance and $\| \mbox{\textbf{P}}^k
\left( (\vec{x}_0,e), \cdot \right) - \mbox{\textbf{P}}^{\infty} (\cdot)
\|_2$ in the special case when $\widehat{Q}$ is defined by (\ref{1.1}).
Analogous results could be established using bounds for random walks
generated by random $m$-cycles.  When $P$ is the uniform distribution on
$G$, the result reduces to Theorem~\ref{1.3}.

\begin{corollary} \label{2.8}
Let \emph{\textbf{W}} be the Markov chain on the complete monomial group
$G~\wr~S_n$ as in Theorem~{\rm \ref{2.3}}, where $\widehat{Q}$ is the probability 
measure on $\widehat{S}_n$ defined at {\rm (\ref{1.1})}.  Let $k = \frac{1}{2} 
n \log n + \frac{1}{4} n \log \left( \frac{1}{p_{\min}} - 1 \right) + \frac{1}{2}
cn$.  Then there exists a universal constant $b > 0$ such that
\[ 
\| \emph{\textbf{P}}^{k}\left( (\vec{x}_0;e), \cdot \right) -
\emph{\textbf{P}}^{\infty}\left( \cdot \right) \|_{\mbox{\rm \scriptsize
TV}} \ \ \leq \ \ \displaystyle \mbox{$\frac{1}{2}$} \ \|
\emph{\textbf{P}}^k \left( (\vec{x}_0,e), \cdot \right) -
\emph{\textbf{P}}^{\infty} (\cdot) \|_2 \ \ \leq \ \ be^{-c} \ \ \
\emph{for all}\ c > 0. 
\]
\end{corollary}

\proof{Proof} Let $\widehat{Q}$ be defined by (\ref{1.1}).  For any set
$J$ with $|J| = j$, it is clear that we have
\[ 
\mu_n(J) \ = \ \left( j/n \right)^2 \ \ \ \mbox{and \ \ \ $d_k(J) \ = \
j! \| Q_{S_j}^{*k} - U_{S_j} \|_2^2$}, 
\]
\noindent
where $Q_{S_j}$ is the measure on $S_j$ induced
by~(\ref{1.1})
and $U_{S_j}$ 
is the uniform distribution on $S_j$.

It then follows from Theorem~\ref{1.2} that there exists a universal
constant $a > 0$ such that $D_k(j) \leq 4a^2 e^{-2c}$ for each $1 \leq j
\leq n$, when $k \geq \frac{1}{2} j \log j + \frac{1}{2} cj$.  Since $n
\geq j$ and $p_{\min} \leq 1/2$, this is also true when $k = \frac{1}{2} n
\log n + \frac{1}{4} n \log \left( \frac{1}{p_{\min}} - 1 \right) +
\frac{1}{2} cn$.  

It then follows from Corollary~\ref{2.5}, with $m = 2$, that
\[ \begin{array}{rcl}

\displaystyle \| \mbox{\textbf{P}}^{k}\left( (\vec{x}_0;e), \cdot \right)
- \mbox{\textbf{P}}^{\infty}\left( \cdot \right) \|_{\mbox{\rm \scriptsize
TV}}^2  &  \leq  &  \displaystyle \mbox{$\frac{1}{4}$} \ \|
\mbox{\textbf{P}}^k \left( (\vec{x}_0,e), \cdot \right) -
\mbox{\textbf{P}}^{\infty} (\cdot) \|_2^2 \vspace{1pc} \\

&  \leq  &  4a^2 e^{-2c} \ + \ e^{-2c} \ \ = \ \ \left( 4a^2 + 1 \right)
e^{-2c}, 

\end{array} \]
\noindent
from which the desired result follows.  \qed

Corollary~\ref{2.8} shows that $k = \frac{1}{2} n \log n + \frac{1}{4} n
\log \left( \mbox{$\frac{1}{p_{\min}}$} - 1 \right) + \frac{1}{2} cn$
steps are sufficient for the $L^2$ distance, and hence also the total
variation distance, to become small.  A lower bound in the $L^2$ distance
can also be derived by examining $n^2 \left( \mbox{$\frac{1}{p_{\min}}$} -
1 \right) \left( 1 - \frac{1}{n} \right)^{4k}$, which is the contribution,
when $j = n-1$ and $m=2$, to the second summation of (\ref{2.6}) from the
proof of Corollary~\ref{2.5}.  In the present context, the second
summation of (\ref{2.6}) is the second summation in the statement of
Theorem~\ref{2.3} with $\mu_n(J) = (|J|/n)^2$.  Notice that $k =
\frac{1}{2} n \log n + \frac{1}{4} n \log \left( 
\mbox{$\frac{1}{p_{\min}}$} - 1 \right) - \frac{1}{2} cn$ steps are
necessary for just this term to become small.

\section{Markov Chains on $(G~\wr~S_n) / (S_r \times S_{n-r})$.} \label{sec3}

We now analyze a very general Markov chain on the homogeneous space
$(G~\wr~S_n) / (S_r \times S_{n-r})$.  It should be noted that, in the
results which follow, there is no essential use of the group structure on
$G$.  So the results of this section extend simply; in general, the Markov
chain of interest is on the set $G^n \times \left( S_n / (S_r \times
S_{n-r}) \right)$.

\subsection{A Class of Chains on $(G~\wr~S_n) / (S_r \times S_{n-r})$.} \label{sec3.1}

Let $[n] := \{1, 2, \ldots, n \}$ and let $[r] := \{1, 2, \ldots, r \}$
where $1 \leq r \leq n/2$.  Recall that the homogeneous space $X = S_n/
(S_r \times S_{n-r})$ can be identified with the set of all ${n \choose
r}$ subsets of size $r$ from $[n]$.  Suppose that $x = \{ i_1, i_2,
\ldots, i_r \} \subseteq [n]$ is such a subset and that $[n] \setminus x =
\{ j_{r+1}, j_{r+2}, \ldots, j_n \}$.  Let $\{ i_{(1)}, i_{(2)}, \ldots,
i_{(k)} \} \subseteq x$ and $\{ j_{(r+1)}, j_{(r+2)}, \ldots, j_{(r+k)} \}
\subseteq [n] \setminus x$ be the sets with all indices, listed in
increasing order, such that $r+1 \leq i_{(\ell)} \leq n$ and $1 \leq
j_{(\ell)} \leq r$ for $1 \leq \ell \leq k$; in the Bernoulli--Laplace
framework, these are the labels of the balls that are no longer in their
respective initial racks.  (Notice that if all the balls are on their
initial racks, then both of these sets are empty.)  To each element $x \in
X$, we can thus correspond a unique permutation
\[ 
(j_{(r+1)}\ i_{(1)}) (j_{(r+2)}\ i_{(2)}) \cdots (j_{(r+k)}\ i_{(k)}) 
\]
\noindent
in $S_n$, which is the product of $k$ (disjoint) transpositions; when this
permutation serves to represent an element of the homogeneous space $X$,
we denote it by $\tilde{\pi}$.  For example, if $x = \{ 2, 4, 8 \} \in X =
S_8/ (S_3 \times S_5)$, then $\tilde{\pi} = (1\ 4) (3\ 8)$.  (If all of
the balls are on their initial racks, then $\tilde{\pi} = e$.)  Notice
that any given $\pi \in S_n$ corresponds to a unique $\tilde{\pi} \in X$;
denote the mapping $\pi \mapsto \tilde{\pi}$ by $R$.  For example, let
$\pi$ be the permutation that sends $(1,2,3,4,5,6,7,8)$ to
$(8,2,4,6,7,1,5,3)$; then $x = \{8,2,4\} = \{2,4,8\}$ and $\tilde{\pi} =
R(\pi) = (1\ 4) (3\ 8)$.

We now modify the concept of augmented permutation introduced in
Section~\ref{sec2.1}.  Rather than the ordered pair of a permutation $\pi \in
S_n$ and a subset $J$ of $F(\pi)$, we now take an augmented permutation to
be the ordered pair of a permutation $\pi \in S_n$ and a subset $J$ of
$F(R(\pi))$.  $\bigg[$ In the above example, $F(R(\pi)) = F(\tilde{\pi}) =
\{2,5,6,7\} \bigg]$.  The necessity of this subtle difference will become
apparent when defining $\widehat{Q}$.  For $\hat{\pi} = (\pi,J) \in
\widehat{S}_n$ (defined in Section~\ref{sec2.1}), define
\[ 
\widetilde{I}(\hat{\pi}) \ \ := \ \ I(R(\pi),J) \ \ = \ \
I(R(T(\hat{\pi})),J). 
\]
\noindent
Thus $\widetilde{I}(\hat{\pi})$ is the union of the set of indices
deranged by $R(T(\hat{\pi}))$ and the subset $J$ of the fixed points of
$R(T(\hat{\pi}))$.

Let~$\widehat{Q}$ be a probability measure on the augmented permutations
$\widehat{S}_n$ satisfying the augmented symmetry property (\ref{2.0}).
Let~$Q$ be as described in Section~\ref{sec2.1}.

Let $\hat{\xi}_1, \hat{\xi}_2, \ldots$ be a sequence of independent
augmented permutations each distributed according to $\widehat{Q}$.  These
correspond uniquely to a sequence $\xi_1, \xi_2, \ldots$ of permutations
each distributed according to $Q$.  Define \textbf{Y} $:= (Y_0, Y_1, Y_2,
\ldots)$ to be the Markov chain on $S_n/ (S_r \times S_{n-r})$ such that
$Y_0 := \tilde{e}$ and $Y_k := R\left( \xi_k Y_{k-1} \right)$ for all $k
\geq 1$.

Let $P$ be a probability measure defined on a finite group $G$ and let
$p_i$ for $1~\leq~i~\leq~|G|$ and $p_{\min} > 0$ be defined as in
Section~\ref{sec2.1}.  Define \textbf{X} $:= (X_0, X_1, X_2, \ldots)$ to be 
the Markov chain on $G^n$ such that $X_0 := \vec{x}_0 = (\chi_1, \ldots,
\chi_n)$ with $\chi_i \in G$ for $1 \leq i \leq n$ and, at each step $k$
for $k \geq 1$, the entries of $X_{k-1}$ whose positions are included in
$I(\hat{\xi}_k)$ are independently changed to an element of $G$
distributed according to $P$.  

Define \textbf{W} $:= (W_0, W_1, W_2, \ldots)$ to be the Markov chain on
$(G~\wr~S_n) / (S_r \times S_{n-r})$ such that $W_k := (X_k; Y_k)$ for all
$k \geq 0$.  Notice that the signed generalization of the classical
Bernoulli--Laplace diffusion model analyzed in Theorem~\ref{1.6}
is a special case of \textbf{W}, with
$P$ being the uniform distribution on $\mathbb{Z}_2$ and $\widehat{Q}$
being defined as at (\ref{1.4}). 

Let \textbf{P}$(\cdot, \cdot)$ be the transition matrix for \textbf{W} and
let \textbf{P}$^{\infty}(\cdot)$ be the stationary distribution for
\textbf{W}.  Notice that \vspace{-0.15in} 
\[ 
\mbox{\textbf{P}}^{\infty} \left( \vec{x};\tilde{\pi} \right) \ = \
\frac{1}{{n \choose r}} \prod_{i=1}^n p_{x_i} 
\]
\noindent
for any $(\vec{x};\tilde{\pi}) \in (G~\wr~S_n) / (S_r \times S_{n-r})$ and
that 
\[ 
\mbox{\textbf{P}} \left( (\vec{x};\tilde{\pi}), 
(\vec{y};\tilde{\sigma}) \right) \ = \ \sum_{\hat{\rho} \in \widehat{S}_n :
R(T(\hat{\rho})\tilde{\pi}) = \tilde{\sigma}} \widehat{Q}(\hat{\rho})
\left[ \prod_{j \in I(\hat{\rho})} p_{y_j} \right] \cdot \left[
\prod_{\ell \not\in I(\hat{\rho})} \mathbb{I}(x_{\ell} = y_{\ell}) \right]
\]
\noindent 
for any $(\vec{x};\tilde{\pi}), (\vec{y};\tilde{\sigma}) \in (G~\wr~S_n) /
(S_r \times S_{n-r})$.  Thus, using the augmented symmetry of
$\widehat{Q}$,
\[ \begin{array}{l}

\displaystyle \mbox{\textbf{P}}^{\infty}\left( \vec{x};\tilde{\pi} \right)
\mbox{\textbf{P}} \left( (\vec{x};\tilde{\pi}), (\vec{y};\tilde{\sigma})
\right) \vspace{1pc} \\

\ \ \ = \displaystyle \left[ \frac{1}{{n \choose r}} \prod_{i=1}^n p_{x_i}
\right] \sum_{\hat{\rho} \in \widehat{S}_n : R(T(\hat{\rho})\tilde{\pi}) =
\tilde{\sigma}} \widehat{Q}(\hat{\rho}) \left[ \prod_{j \in I(\hat{\rho})}
p_{y_j} \right] \cdot \left[ \prod_{\ell \not\in I(\hat{\rho})}
\mathbb{I}(x_{\ell} = y_{\ell}) \right] \vspace{1pc} \\

\ \ \ = \displaystyle \sum_{\hat{\rho} \in \widehat{S}_n :
R(T(\hat{\rho})\tilde{\pi}) = \tilde{\sigma}} \widehat{Q}(\hat{\rho})
\left[ \frac{1}{{n \choose r}} \prod_{i \in I(\hat{\rho})} p_{x_i}
\prod_{i \not\in I(\hat{\rho})} p_{x_i} \right] \cdot \left[ \prod_{j \in
I(\hat{\rho})} p_{y_j} \right] \cdot \left[ \prod_{\ell \not\in 
I(\hat{\rho})} \mathbb{I}(x_{\ell} = y_{\ell}) \right] \vspace{1pc} \\

\ \ \ = \displaystyle \sum_{\hat{\rho} \in \widehat{S}_n :
R(T(\hat{\rho})\tilde{\sigma}) = \tilde{\pi}} \widehat{Q}(\hat{\rho}) 
\left[ \frac{1}{{n \choose r}} \prod_{i \in I(\hat{\rho})} p_{x_i} 
\prod_{j \not\in I(\hat{\rho})} p_{y_j} \right] \cdot \left[ \prod_{j \in
I(\hat{\rho})} p_{y_j} \right] \cdot \left[ \prod_{\ell \not\in
I(\hat{\rho})} \mathbb{I}(y_{\ell} = x_{\ell}) \right] \vspace{1pc} \\

\ \ \ = \displaystyle \left[ \frac{1}{{n \choose r}} \prod_{j=1}^n p_{y_j}
\right] \sum_{\hat{\rho} \in \widehat{S}_n : 
R(T(\hat{\rho})\tilde{\sigma}) = \tilde{\pi}} \widehat{Q}(\hat{\rho})
\left[ \prod_{i \in I(\hat{\rho})} p_{x_i} \right] \cdot \left[
\prod_{\ell \not\in I(\hat{\rho})} \mathbb{I}(y_{\ell} = x_{\ell}) \right]
\vspace{1pc} \\

\ \ \ = \displaystyle \mbox{\textbf{P}}^{\infty}\left( 
\vec{y};\tilde{\sigma} \right) \mbox{\textbf{P}} \left( 
(\vec{y};\tilde{\sigma}), (\vec{x};\tilde{\pi}) \right).

\end{array} \]
\noindent
Therefore, \textbf{P} is reversible, which is a necessary condition in
order to apply the comparison technique of Diaconis and Saloff-Coste (1993b).

\subsection{Convergence to Stationarity:\ Main Result.} \label{sec3.2}

For any $J \subseteq [n]$, let $X^{(J)}$ be the homogeneous space
$S_{(J)}/ \left( S_{(J \cap [r])} \times S_{(J \cap ([n] \setminus
[r]))} \right)$, where $S_{(J')}$ is the subgroup of $S_n$ consisting of
those $\sigma \in S_n$ with $[n] \setminus F(\sigma) \subseteq J'$.
As in Section~\ref{sec3.1},
let~$\widehat{Q}$ be a probability measure on the augmented permutations
$\widehat{S}_n$ satisfying the augmented symmetry property (\ref{2.0}).

Let $Q$ and $Q_{S_{(J)}}$ be as described in Sections~\ref{sec2.1} and
\ref{sec2.2}.  For notational purposes, let 
\num \begin{equation} \label{3.1}
\tilde{\mu}_n(J) \ \ := \ \ \widehat{Q}\{ \hat{\sigma} \in \widehat{S}_n :
\widetilde{I}(\hat{\sigma}) \subseteq J\}.
\end{equation}

Let $\widetilde{Q}_{X^{(J)}}$ be the probability measure on $X^{(J)}$
induced (as described in Section 2.2 of Schoolfield (1999b)) by
$Q_{S_{(J)}}$.  Also let $\widetilde{U}_{X^{(J)}}$ be the uniform measure
on $X^{(J)}$.  For notational purposes, let 
\num \begin{equation} \label{3.2}
\tilde{d}_k(J) \ \ := \ \ \mbox{${|J| \choose |J \cap [r]|}$} \|
\widetilde{Q^{*k}}_{X^{(J)}} - \widetilde{U}_{X^{(J)}} \|_2^2.
\end{equation}
\bigskip

\emph{Example.\/}\ Let $\widehat{Q}$ be defined as at (\ref{1.4}).  Then $\widehat{Q}$ 
satisfies the augmented symmetry property (\ref{2.0}).  In the
Bernoulli--Laplace framework, the elements $\widehat{Q}(\kappa,\{j\})$
and $\widehat{Q}(\kappa,\{i,j\})$ leave the balls on their current racks,
but single out one or two of them, respectively; the element
$\widehat{Q}(\tau \kappa,\emptyset)$ switches two balls between the racks.
In Corollary~\ref{3.8} we will be using~$\widehat{Q}$ to define a Markov
chain on $(G~\wr~S_n) / (S_r \times S_{n-r})$
which is a generalization of the Markov chain analyzed in Theorem~\ref{1.6}.

It is also easy to verify that $\widehat{Q}_{S_{(J)}}$ is the probability 
measure defined at (\ref{1.4}), but with $[r]$ and $[n] \setminus [r]$
changed to $J \cap [r]$ and $J \cap ([n] \setminus [r])$, respectively.
Thus, roughly put, our generalization of the Markov chain analyzed in
Theorem~\ref{1.6}, conditionally restricted to the indices in $J$, gives
a Markov chain on $\left( G \wr S_{(J)} \right) / \left( S_{(J \cap [r])}
\times S_{(J \cap ([n] \setminus [r]))} \right)$ ``as if $J$ were the only
indices.''
\bigskip

The following result establishes an upper bound on the total variation
distance by deriving an exact formula for $\| \mbox{\textbf{P}}^{k}\left(
(\vec{x}_0;\tilde{e}), \cdot \right) - \mbox{\textbf{P}}^{\infty} (\cdot)
\|_2^2$.

\begin{theorem} \label{3.3}
Let \emph{\textbf{W}} be the Markov chain on the homogeneous space \linebreak
$(G~\wr~S_n) / (S_r \times S_{n-r})$ defined in Section~{\rm \ref{sec3.1}}.  Then
\[ \begin{array}{rcl}

\displaystyle \| \emph{\textbf{P}}^{k}\left( (\vec{x}_0;\tilde{e}), \cdot
\right) - \emph{\textbf{P}}^{\infty}\left( \cdot \right) \|_{\mbox{\rm
\scriptsize TV}}^2  &  \leq  &  \displaystyle \mbox{$\frac{1}{4}$} \ \|
\emph{\textbf{P}}^{k}\left( (\vec{x}_0;\tilde{e}), \cdot \right) -
\emph{\textbf{P}}^{\infty} (\cdot) \|_2^2 \vspace{1pc} \\

&  =  &  \displaystyle \mbox{$\frac{1}{4}$} \sum_{J: J \subseteq [n]}
\frac{{n \choose r}}{{|J| \choose |J \cap [r]|}} \left[ \prod_{i \not\in
J} \left( \mbox{$\frac{1}{p_{\chi_i}}$} - 1 \right) \right] 
\tilde{\mu}_n(J)^{2k} \ \tilde{d}_k(J) \vspace{1pc} \\

&     &  \displaystyle + \ \ \mbox{$\frac{1}{4}$} \sum_{J: J \subsetneq
[n]} \frac{{n \choose r}}{{|J| \choose |J \cap [r]|}} \left[ \prod_{i
\not\in J} \left( \mbox{$\frac{1}{p_{\chi_i}}$} - 1 \right) \right]
\tilde{\mu}_n(J)^{2k},

\end{array} \]
\noindent
where $\tilde{\mu}_n(J)$ and $\tilde{d}_k(J)$ are defined at {\rm (\ref{3.1})}
and {\rm (\ref{3.2})}, respectively.
\end{theorem}

\proof{Proof} For each $k \geq 1$, let $\displaystyle H_k := \bigcup_{\ell
= 1}^k \widetilde{I}(\hat{\xi}_{\ell}) \subseteq [n]$.  For any given $w =
(\vec{x};\tilde{\pi}) \in (G~\wr~S_n) / (S_r \times S_{n-r})$, let $A
\subseteq [n]$ be the set of indices such that $x_i \neq \chi_i$, where
$x_i$ is the $i$th entry of $\vec{x}$ and $\chi_i$ is the $i$th entry of
$\vec{x}_0$, and let $B = [n] \setminus F(\tilde{\pi})$ be the set of
indices deranged by $\tilde{\pi}$.  Notice that $H_k \supseteq A \cup B$.

The proof continues exactly as in the proof of Theorem~\ref{2.3} to
determine that
\[ 
\displaystyle \mathbb{P}\left( W_k = (\vec{x};\tilde{\pi}) \right) \ =
\ \sum_{J: B \subseteq J \subseteq [n]} \left( -1 \right)^{|J|}
\tilde{\mu}_n(J)^{k} \ \mathbb{P}\left( Y_k = \tilde{\pi} \ |\ H_k
\subseteq J \right) \prod_{i=1}^n \left[ 1 - \mathbb{I}_{A \cup J}(i) -
p_{x_i} \right]. 
\]

In particular, when $(\vec{x};\tilde{\pi}) = (\vec{x}_0;\tilde{e})$, we
have $A = \emptyset = B$ and
\[ \begin{array}{rcl}

\displaystyle \mathbb{P}\left( W_k = (\vec{x}_0;\tilde{e}) \right)  &  =
&  \displaystyle \sum_{J: J \subseteq [n]} \left( -1 \right)^{|J|}
\tilde{\mu}_n(J)^k \ \mathbb{P}\left( Y_k = \tilde{e} \ |\ H_k \subseteq J
\right) \prod_{i=1}^n \left[ 1 - \mathbb{I}_J(i) - p_{\chi_i} \right]
\vspace{1pc} \\

&  =  &  \displaystyle \left[ \prod_{i=1}^n p_{\chi_i} \right] \sum_{J: J
\subseteq [n]} \tilde{\mu}_n(J)^k \ \mathbb{P}\left( Y_k = \tilde{e} \ |\
H_k \subseteq J \right) \prod_{i \not\in J} \left(
\mbox{$\frac{1}{p_{\chi_i}}$} - 1 \right). 

\end{array} \]
\noindent
Notice that $\displaystyle \{ H_k \subseteq J \} = \bigcap_{\ell=1}^k
\left\{ \widetilde{I}(\hat{\xi}_{\ell}) \subseteq J \right\}$ for any $k$
and $J$.  So $\mathcal{L}$ $\left( (Y_0, Y_1, \ldots, Y_k \ | \ H_k
\subseteq J) \right)$ is the law of a Markov chain on $S_n/ (S_r \times
S_{n-r})$ (through step $k$) with step distribution $Q_{X^{(J)}}$.  Thus,
using the reversibility of \textbf{P} and the symmetry of $Q_{X^{(J)}}$,
\[ \begin{array}{l}

\| \mbox{\textbf{P}}^{k}\left( (\vec{x}_0;\tilde{e}), \cdot \right) -
\mbox{\textbf{P}}^{\infty} (\cdot) \|_2^2 \ \ = \ \ \displaystyle \frac{{n
\choose r}}{\prod_{i=1}^n p_{\chi_i}} \mbox{\textbf{P}}^{2k}\left(
(\vec{x}_0;\tilde{e}), (\vec{x}_0;\tilde{e}) \right) \ - \ 1 \vspace{1pc}
\\

\ \ \ = \displaystyle \ {n \choose r} \sum_{J: J \subseteq [n]} \ \left[
\prod_{i \not\in J} \left( \mbox{$\frac{1}{p_{\chi_i}}$} - 1 \right)
\right] \ \tilde{\mu}_n(J)^{2k} \ \mathbb{P}\left( Y_{2k} = \tilde{e} \ |\
H_{2k} \subseteq J \right) \ \ - \ \ 1 \vspace{1pc} \\

\ \ \ = \displaystyle \ {n \choose r} \sum_{J: J \subseteq [n]} \ \left[
\prod_{i \not\in J} \left( \mbox{$\frac{1}{p_{\chi_i}}$} - 1 \right)
\right] \ \tilde{\mu}_n(J)^{2k} \ \left[ \| \widetilde{Q^{*k}}_{X^{(J)}} -
\widetilde{U}_{X^{(J)}} \|_2^2 + \frac{1}{{|J| \choose |J \cap [r]|}}
\right] \ \ - \ \ 1 \vspace{1pc} \\

\ \ \ = \displaystyle \ {n \choose r} \sum_{J: J \subseteq [n]} \ \left[
\prod_{i \not\in J} \left( \mbox{$\frac{1}{p_{\chi_i}}$} - 1 \right)
\right] \ \tilde{\mu}_n(J)^{2k} \ \frac{1}{{|J| \choose |J \cap [r]|}} \
\left( \tilde{d}_k(J) + 1 \right) \ \ - \ \ 1 \vspace{1pc} \\

\ \ \ = \displaystyle \ \sum_{J: J \subseteq [n]} \frac{{n \choose 
r}}{{|J| \choose |J \cap [r]|}} \left[ \prod_{i \not\in J} \left( 
\mbox{$\frac{1}{p_{\chi_i}}$} - 1 \right) \right] \tilde{\mu}_n(J)^{2k} \
\tilde{d}_k(J) \vspace{1pc} \\

\ \ \ \ \ \ \ \ + \displaystyle \ \sum_{J: J \subsetneq [n]} \frac{{n
\choose r}}{{|J| \choose |J \cap [r]|}} \left[ \prod_{i \not\in J} \left(
\mbox{$\frac{1}{p_{\chi_i}}$} - 1 \right) \right] \tilde{\mu}_n(J)^{2k},

\end{array} \]
\noindent
from which the desired result follows. $\qed$

\subsection{Corollaries.} \label{sec3.3}

We now establish several corollaries to our main result.

\begin{corollary} \label{3.4}
Let \emph{\textbf{W}} be the Markov chain on the homogeneous space \linebreak
$(G~\wr~S_n) / (S_r~\times~S_{n-r})$ as in Theorem~{\rm \ref{3.3}}.  For $0 \leq
j \leq n$, let
\[ 
\widetilde{M}_n(j) := \max \left\{ \tilde{\mu}_n(J) : |J| = j \right\}
\ \ \ and \ \ \ \widetilde{D}_k(j) := \max \left\{ \tilde{d}_k(J) : |J| =
j \right\}. 
\]
\noindent
Also let
\[ 
\widetilde{B}(n,k) := \max \left\{ \widetilde{D}_k(j) : 0 \leq j \leq n
\right\} = \max \left\{ \tilde{d}_k(J) : J \subseteq [n] \right\}. 
\]
\noindent
Then
\[ \begin{array}{l}

\displaystyle \| \emph{\textbf{P}}^{k}\left( (\vec{x}_0;\tilde{e}), \cdot
\right) - \emph{\textbf{P}}^{\infty}\left( \cdot \right) \|_{\mbox{\rm
\scriptsize TV}}^2 \ \ \leq \ \ \displaystyle \mbox{$\frac{1}{4}$} \ \|
\emph{\textbf{P}}^k \left( (\vec{x}_0,\tilde{e}) - \cdot \right),
\emph{\textbf{P}}^{\infty} (\cdot) \|_2^2 \vspace{1pc} \\

\ \ \ \leq \displaystyle \ \mbox{$\frac{1}{4}$} \ \widetilde{B}(n,k)
\sum_{i=0}^r \sum_{j=0}^{n-r} {r \choose i} {n-r \choose j} \frac{{n
\choose r}}{{i+j \choose i}} \left( \mbox{$\frac{1}{p_{\min}}$} - 1
\right)^{n-(i+j)} \widetilde{M}_n(j)^{2k} \vspace{1pc} \\

\ \ \ \ \ \ \ \ + \displaystyle \ \mbox{$\frac{1}{4}$} \sum_{i=0}^r
\sum_{j=0}^{n-r} {r \choose i} {n-r \choose j} \frac{{n \choose
r}}{{i+j \choose i}} \left( \mbox{$\frac{1}{p_{\min}}$} - 1
\right)^{n-(i+j)} \widetilde{M}_n(j)^{2k},

\end{array} \]
\noindent
where the last sum must be modified to exclude the term for $i = r$ and $j
= n-r$.
\end{corollary}

\proof{Proof} The proof is analogous to that of Corollary~\ref{2.4}. \qed

\begin{corollary} \label{3.5}
In addition to the assumptions of Theorem~{\rm \ref{3.3}} and
Corollary~{\rm \ref{3.4}}, suppose that there exists $m > 0$ such that
$\displaystyle \widetilde{M}_n(j) \leq \left( j/n \right)^m$ for
\emph{all} $0 \leq j \leq n$.  Let $k \geq \frac{1}{2m} n \left( \log n +
\log \left( \mbox{$\frac{1}{p_{\min}}$} - 1 \right) + c \right)$.  Then
\[ 
\| \emph{\textbf{P}}^{k}\left( (\vec{x}_0;\tilde{e}), \cdot \right),
\emph{\textbf{P}}^{\infty}\left( \cdot \right) \|_{\mbox{\rm \scriptsize
TV}} \ \ \leq \ \ \displaystyle \mbox{$\frac{1}{2}$} \ \|
\emph{\textbf{P}}^k \left( (\vec{x}_0,\tilde{e}), \cdot \right) -
\emph{\textbf{P}}^{\infty} (\cdot) \|_2 \ \ \leq \ \ 2 \left(
\widetilde{B}\left(n,k\right) \ + \ e^{-c} \right)^{1/2}. 
\]
\end{corollary}

\proof{Proof} It follows from Corollary~\ref{3.4} that
\num \begin{equation} \label{3.6}
\begin{array}{l}

\displaystyle \| \mbox{\textbf{P}}^{k}\left( (\vec{x}_0;\tilde{e}), \cdot
\right) - \mbox{\textbf{P}}^{\infty}\left( \cdot \right) \|_{\mbox{\rm
\scriptsize TV}}^2 \ \ \leq \ \ \displaystyle \mbox{$\frac{1}{4}$} \ \|
\mbox{\textbf{P}}^{k}\left( (\vec{x}_0;\tilde{e}), \cdot \right) -
\mbox{\textbf{P}}^{\infty} (\cdot) \|_2^2 \nonumber \vspace{1pc} \\

\ \ \ \leq \displaystyle \ \mbox{$\frac{1}{4}$} \
\widetilde{B}\left(n,k\right) \sum_{i=0}^r \sum_{j=0}^{n-r} {r \choose i}
{n-r \choose j} \frac{{n \choose r}}{{i+j \choose i}} \left(
\mbox{$\frac{1}{p_{\min}}$} - 1 \right)^{n-(i+j)} \left( \frac{i+j}{n} 
\right)^{2km} \nonumber \vspace{1pc} \\

\ \ \ \ \ \ \ \ + \displaystyle \ \mbox{$\frac{1}{4}$} \sum_{i=0}^r
\sum_{j=0}^{n-r} {r \choose i} {n-r \choose j} \frac{{n \choose
r}}{{i+j \choose i}} \left( \mbox{$\frac{1}{p_{\min}}$} - 1
\right)^{n-(i+j)} \left( \frac{i+j}{n} \right)^{2km},

\end{array}
\end{equation}
\noindent
where the last sum must be modified to exclude the term for $i = r$ and $j
= n-r$.  Notice that
\[ 
{r \choose i} {n-r \choose j} \frac{{n \choose r}}{{i+j \choose i}} \ =
\ {n \choose i+j} { n-(i+j) \choose r-i}. 
\]
\noindent
Thus if we put $j' = i+j$ and change the order of summation we have
(enacting now the required modification)
\[ \begin{array}{l}

\displaystyle \| \mbox{\textbf{P}}^{k}\left( (\vec{x}_0;\tilde{e}), \cdot
\right), \mbox{\textbf{P}}^{\infty}\left( \cdot \right) \|_{\mbox{\rm
\scriptsize TV}}^2 \ \ \leq \ \ \displaystyle \mbox{$\frac{1}{4}$} \ \|
\mbox{\textbf{P}}^{k}\left( (\vec{x}_0;\tilde{e}), \cdot \right) -
\mbox{\textbf{P}}^{\infty} (\cdot) \|_2^2 \vspace{1pc} \\

\ \ \ \leq \displaystyle \ \mbox{$\frac{1}{4}$} \
\widetilde{B}\left(n,k\right) \sum_{j=0}^n {n \choose j} \left(
\mbox{$\frac{1}{p_{\min}}$} - 1 \right)^{n-j} \left( \frac{j}{n}
\right)^{2km} \sum_{i=\ell \vee (r-(n-j))}^{r \wedge (j-\ell)} {n-j 
\choose r-i} \vspace{1pc} \\

\ \ \ \ \ \ \ \ + \displaystyle \ \mbox{$\frac{1}{4}$} \sum_{j=0}^{n-1} {n
\choose j} \left( \mbox{$\frac{1}{p_{\min}}$} - 1 \right)^{n-j} \left(
\frac{j}{n} \right)^{2km} \sum_{i=\ell \vee (r-(n-j))}^{r \wedge (j-\ell)}
{n-j \choose r-i}.

\end{array} \]

Of course $\displaystyle \sum_{i=\ell \vee (r-(n-j))}^{r \wedge (j-\ell)}
\mbox{${n-j \choose r-i}$} \ \leq \ 2^{n-j}$.  If we then let $i = n-j$,
the upper bound becomes
\[ \begin{array}{l}

\displaystyle \| \mbox{\textbf{P}}^{k}\left( (\vec{x}_0;\tilde{e}), \cdot
\right) - \mbox{\textbf{P}}^{\infty}\left( \cdot \right) \|_{\mbox{\rm
\scriptsize TV}}^2 \ \ \leq \ \ \displaystyle \mbox{$\frac{1}{4}$} \ \|
\mbox{\textbf{P}}^{k}\left( (\vec{x}_0;\tilde{e}), \cdot \right) -
\mbox{\textbf{P}}^{\infty} (\cdot) \|_2^2 \vspace{1pc} \\

\ \ \ \leq \displaystyle \ \mbox{$\frac{1}{4}$} \
\widetilde{B}\left(n,k\right) \sum_{i=0}^n 2^i {n \choose i} \left( 
\mbox{$\frac{1}{p_{\min}}$} - 1 \right)^i \left( 1 - \mbox{$\frac{i}{n}$}
\right)^{2km} \vspace{1pc} \\

\ \ \ \ \ \ \ \ + \displaystyle \ \mbox{$\frac{1}{4}$} \sum_{i=1}^n 2^i {n
\choose i} \left( \mbox{$\frac{1}{p_{\min}}$} - 1 \right)^i \left( 1 -
\mbox{$\frac{i}{n}$} \right)^{2km} \vspace{1pc} \\

\ \ \ \leq  \displaystyle \ \mbox{$\frac{1}{4}$} \ 
\widetilde{B}\left(n,k\right) \sum_{i=0}^n \frac{1}{i!} (2n)^i \left(
\mbox{$\frac{1}{p_{\min}}$} - 1 \right)^i e^{-2ikm/n} \ \ + \ \
\mbox{$\frac{1}{4}$} \sum_{i=1}^n \frac{1}{i!} (2n)^i \left(
\mbox{$\frac{1}{p_{\min}}$} - 1 \right)^i e^{-2ikm/n}. 

\end{array} \]
\noindent
Notice that if $k \geq \frac{1}{2m} n \left( \log n + \log \left(
\mbox{$\frac{1}{p_{\min}}$} - 1 \right) + c \right)$, then
\[ 
e^{-2ikm/n} \ \leq \ \left[ \frac{e^{-c}}{\left(
\mbox{$\frac{1}{p_{\min}}$} - 1 \right) n} \right]^i, 
\]
\noindent
from which it follows that
\[ \begin{array}{rcl}

\displaystyle \| \mbox{\textbf{P}}^{k}\left( (\vec{x}_0;\tilde{e}), \cdot
\right) - \mbox{\textbf{P}}^{\infty}\left( \cdot \right) \|_{\mbox{\rm
\scriptsize TV}}^2  &  \leq  &  \displaystyle \mbox{$\frac{1}{4}$} \ \|
\mbox{\textbf{P}}^{k}\left( (\vec{x}_0;\tilde{e}), \cdot \right) -
\mbox{\textbf{P}}^{\infty} (\cdot) \|_2^2 \vspace{1pc} \\

&  \leq  &  \displaystyle \mbox{$\frac{1}{4}$} \
\widetilde{B}\left(n,k\right) \sum_{i=0}^n \frac{1}{i!} \left( 2e^{-c}
\right)^i \ \ + \ \ \mbox{$\frac{1}{4}$} \sum_{i=1}^n \frac{1}{i!} \left(
2e^{-c} \right)^i \vspace{1pc} \\

&  \leq  &  \displaystyle \mbox{$\frac{1}{4}$} \
\widetilde{B}\left(n,k\right) \exp\left( 2e^{-c} \right) \ \ + \ \ 
\mbox{$\frac{1}{2}$} e^{-c} \exp\left( 2e^{-c} \right).

\end{array} \]
\noindent
Since $c > 0$, we have $\exp\left(2e^{-c}\right) < e^2$.  Therefore
\[ 
\displaystyle \| \mbox{\textbf{P}}^{k}\left( (\vec{x}_0;\tilde{e}),
\cdot \right) - \mbox{\textbf{P}}^{\infty}\left( \cdot \right)
\|_{\mbox{\rm \scriptsize TV}}^2 \ \ \leq \ \ \displaystyle
\mbox{$\frac{1}{4}$} \ \| \mbox{\textbf{P}}^{k}\left( 
(\vec{x}_0;\tilde{e}), \cdot \right) - \mbox{\textbf{P}}^{\infty} (\cdot)
\|_2^2 \ \ \leq \ \ 4 \left( \widetilde{B}\left(n,k\right) \ + \ e^{-c}
\right), 
\]
\noindent
from which the desired result follows.  $\qed$

\begin{corollary} \label{3.7}
In addition to the assumptions of Theorem~{\rm \ref{3.3}} and
Corollary~{\rm \ref{3.4}}, suppose that a set with the distribution of
$\widetilde{I}(\hat{\sigma})$ when $\hat{\sigma}$ has distribution
$\widehat{Q}$ can be constructed by first choosing a set size $0 < \ell
\leq n$ according to a probability mass function $f_n( \cdot )$ and then
choosing a set $L$ with $|L| = \ell$ \emph{uniformly} among all such
choices.  Let $k \geq \frac{1}{2} n \left( \log n + \log \left(
\mbox{$\frac{1}{p_{\min}}$} - 1 \right) + c \right)$. Then 
\[ 
\| \emph{\textbf{P}}^{k}\left( (\vec{x}_0;\tilde{e}), \cdot \right) -
\emph{\textbf{P}}^{\infty}\left( \cdot \right) \|_{\mbox{\rm \scriptsize
TV}} \ \ \leq \ \ \displaystyle \mbox{$\frac{1}{2}$} \ \|
\emph{\textbf{P}}^k \left( (\vec{x}_0,\tilde{e}), \cdot \right) -
\emph{\textbf{P}}^{\infty} (\cdot) \|_2 \ \ \leq \ \ 2 \left(
\widetilde{B}(n,k) \ + \ e^{-c} \right)^{1/2}. 
\]
\end{corollary}

\proof{Proof} The proof is analogous to that of Corollary~\ref{2.7}. \qed

Theorem~\ref{3.3}, and its subsequent corollaries, can be used to bound
the distance to stationarity of many different Markov chains \textbf{W} on
$(G~\wr~S_n) / (S_r \times S_{n-r})$ for which bounds on the $L^2$
distance to uniformity for the related Markov chains on $S_{i+j}/ (S_i
\times S_j)$ for $0 \leq i \leq r$ and $0 \leq j \leq n-r$ are known.  As
an example, the following result establishes an upper bound on both the
total variation distance and $\| \mbox{\textbf{P}}^k \left( 
(\vec{x}_0,\tilde{e}), \cdot \right) - \mbox{\textbf{P}}^{\infty} (\cdot)
\|_2$ in the special case when $\widehat{Q}$ is defined by (\ref{1.4}).
This corollary actually fits the framework of Corollary~\ref{3.7}, but the
result is better than that which would have been determined by merely
applying Corollary~\ref{3.7}.  When $G = \mathbb{Z}_2$ and $P$ is the
uniform distribution on $G$, the result reduces to
Theorem~\ref{1.6}.

\begin{corollary} \label{3.8}
Let \emph{\textbf{W}} be the Markov chain on the homogeneous space \linebreak
$(G~\wr~S_n) / (S_r~\times~S_{n-r})$ as in Theorem~{\rm \ref{3.3}}, where 
$\widehat{Q}$ is the probability measure on $\widehat{S}_n$ defined at 
{\rm (\ref{1.4})}.  Let $k = \frac{1}{4} n \left(\log n + \log \left( 
\frac{1}{p_{\min}} - 1 \right) + c\right)$.  Then there exists a universal 
constant $b > 0$ such that
\[ 
\| \emph{\textbf{P}}^{k}\left( (\vec{x}_0;\tilde{e}), \cdot \right) -
\emph{\textbf{P}}^{\infty}\left( \cdot \right) \|_{\mbox{\rm \scriptsize
TV}} \ \ \leq \ \ \displaystyle \mbox{$\frac{1}{2}$} \ \|
\emph{\textbf{P}}^{k}\left( (\vec{x}_0;\tilde{e}), \cdot \right) -
\emph{\textbf{P}}^{\infty} (\cdot) \|_2 \ \ \leq \ \ be^{-c/2} \ \ \
\emph{for\ all}\ c > 0. 
\]
\end{corollary}

\proof{Proof} The proof is analogous to that of Corollary~\ref{2.8}. \qed

Corollary~\ref{3.8} shows that $k = \frac{1}{4} n \left( \log n + \log
\left( \mbox{$\frac{1}{p_{\min}}$} - 1 \right) + c \right)$ steps are
sufficient for the $L^2$ distance, and hence also the total variation
distance, to become small.  A lower bound in the $L^2$ distance can also
be derived by examining $2n \left( \mbox{$\frac{1}{p_{\min}}$} - 1 \right)
\left( 1 - \frac{1}{n} \right)^{4k}$, which is the contribution, when $i+j
= n-1$ and $m=2$, to the second summation of (\ref{3.6}) from the proof of
Corollary~\ref{3.5}.  In the present context, the second summation of
(\ref{3.6}) is the second summation in the statement of Theorem~\ref{3.3}
with $\tilde{\mu}_n(J) = (|J|/n)^2$.  Notice that $k = \frac{1}{4} n
\left( \log n + \log \left( \mbox{$\frac{1}{p_{\min}}$} - 1 \right) - c
\right)$ steps are necessary for just this term to become small.

\section*{Acknowledgments.} 
This paper
derived from
a portion of the second author's Ph.D. dissertation in
the Department of Mathematical Sciences at the Johns Hopkins University.

\twospace

\Line{
\AOPaddress{James Allen Fill\\
Department of Mathematical Sciences\\
The Johns Hopkins University\\
3400 N. Charles Street\\
Baltimore, Maryland 21218-2682\\
e-mail: {\tt jimfill@jhu.edu}\\
URL: {\tt http://www.mts.jhu.edu/\~{}fill/}} \hfill
\AOPaddress{Clyde H. Schoolfield, Jr.\\
Department of Statistics\\
Harvard University\\
One Oxford Street\\
Cambridge, Massachusetts 02138\\
e-mail: {\tt clyde@stat.harvard.edu}\\
URL: {\tt http://www.fas.harvard.edu/\~{}chschool/}}}

\end{document}